\documentclass{article}
\usepackage{amsthm,amsfonts,amsmath}
\usepackage{amscd,amssymb,latexsym,epsfig}
\usepackage[all]{xy}
\title{Three approaches towards
       Floer homology of cotangent bundles}
\author{Joa Weber\\ \\
        Universit\"at M\"unchen}
\date{28 February 2005\footnote{Last
      revision: 14 January 2006}}
\newtheorem{theorem}{Theorem}[section]

\newtheorem{proposition}[theorem]{Proposition}

\newtheorem{remark}[theorem]{Remark}

\renewcommand{\1}{{{\mathchoice {\rm 1\mskip-4mu l} {\rm 1\mskip-4mu l}
{\rm 1\mskip-4.5mu l} {\rm 1\mskip-5mu l}}}}
\renewcommand{\graph}{{\rm graph }}  
\newcommand{\C}{{\mathbb{C}}}

\newcommand{\N}{{\mathbb{N}}}

\newcommand{\R}{{\mathbb{R}}}

\newcommand{\T}{{\mathbb{T}}}
\newcommand{\Z}{{\mathbb{Z}}}
\newcommand{\Aa}{{\mathcal{A}}}   

\newcommand{\Dd}{{\mathcal{D}}}

\newcommand{\Ff}{{\mathcal{F}}}
\newcommand{\Gg}{{\mathcal{G}}}   

\newcommand{\Ll}{{\mathcal{L}}}   
\newcommand{\Mm}{{\mathcal{M}}}   

\newcommand{\Pp}{{\mathcal{P}}}
\newcommand{\Qq}{{\mathcal{Q}}}

\newcommand{\Ss}{{\mathcal{S}}}
\newcommand{\Tt}{{\mathcal{T}}}

\newcommand{\Vv}{{\mathcal{V}}}

\newcommand{\im}{{\rm im\: }}      
\newcommand{\id}{{\rm id}}         

\newcommand{\IND}{{\rm ind}}       

\newcommand{\grad}{{\rm grad }}    
\newcommand{\Fix}{{\rm Fix}}          
\newcommand{\HF}{{\rm HF}}            
\newcommand{\HI}{{\rm HI}}            
\newcommand{\Ho}{{\rm H}}             
\newcommand{\CM}{{\rm CM}}            
\newcommand{\HM}{{\rm HM}}            
\newcommand{\CF}{{\rm CF}}            

\newcommand{\norm}{{\rm norm}}

\newcommand{\eps}{{\varepsilon}}

\newcommand{\Om}{{\Omega}}

\newcommand{\g}{{\mathfrak g}}     
\newcommand{\Cinf}{C^{\infty}}

\newcommand{\inner}[2]{\langle #1, #2\rangle}

\def\NABLA#1{{\mathop{\nabla\kern-.5ex\lower1ex\hbox{$#1$}}}}
\def\Nabla#1{\nabla\kern-.5ex{}_{#1}}
\def\Tabla#1{\Tilde\nabla\kern-.5ex{}_{#1}}
\def\abs#1{\mathopen|#1\mathclose|}   
\def\Abs#1{\left|#1\right|}            
\def\norm#1{\mathopen\|#1\mathclose\|}
\def\Norm#1{\left\|#1\right\|}

\renewcommand{\Tilde}{\widetilde}

\newcommand{\p}{{\partial}}

\begin{document}

\maketitle


\begin{abstract}
Consider the cotangent bundle
of a closed Riemannian manifold
and an almost complex
structure close
to the one induced by
the Riemannian metric.
For Hamiltonians which
grow for instance
quadratically in the fibers
outside of a compact set,
one can define
Floer homology and show that it is naturally
isomorphic to singular homology
of the free loop space.
We review the three isomorphisms
constructed by Viterbo~\cite{V96},
Salamon-Weber~\cite{JOA3} and
Abbondandolo-Schwarz~\cite{AS04}.
The theory is illustrated
by calculating
Morse and Floer homology
in case of the euclidean 
$n$-torus.
Applications include
existence of noncontractible
periodic orbits
of compactly supported
Hamiltonians on open unit
disc cotangent bundles
which are sufficiently large
over the zero section.
\end{abstract}


\section{Chain group and boundary operators}
\label{sec:intro}

Let $M$ be a closed smooth manifold and
fix a Riemannian metric. Let $\nabla$ be
the associated Levi-Civita connection.
This endows the free loop space
$\Ll M=C^\infty(S^1,M)$
with an $L^2$ and a $W^{1,2}$ metric
given by
$$
     \langle\xi,\eta\rangle_{L^2}
     =\int_0^1\langle\xi(t),\eta(t)\rangle 
     \: dt,\quad
     \langle\xi,\eta\rangle_{W^{1,2}}
     =\langle\xi,\eta\rangle_{L^2} +
     \langle\Nabla{t}\xi,
     \Nabla{t}\eta\rangle_{L^2},
$$
where $\xi$ and $\eta$ are
smooth vector fields
along $x\in\Ll M$.
Here and throughout we identify
$S^1=\R/\Z$ and think of $x\in\Ll M$
as a smooth map $x:\R\to M$
which satisfies $x(t+1)=x(t)$.
Fix a time-dependent function
$V\in C^\infty(S^1\times M)$
and set $V_t(q):=V(t,q)$.
The \emph{classical action functional
on $\Ll M$} is defined by
$$
     \Ss_V(x)
     :=\int_0^1\left(\frac{1}{2}
     \abs{\dot x(t)}^2-V_t(x(t))\right) dt.
$$
The set $\Pp(V)$ of critical points
consists of the $1$-periodic
solutions of the ODE
\begin{equation}\label{eq:crit}
     \Nabla{t} \p_t x=-\nabla V_t(x).
\end{equation}
Here $\nabla V_t$ denotes the gradient.
These solutions are called 
\emph{perturbed closed geodesics}.
Two features make the functional $\Ss_V$
accessible to standard variational methods,
boundedness from below
and finiteness of the
\emph{Morse index}\footnote{The dimension
of the largest subspace
on which the Hessian is negative definite.}
$\IND_V(x)$ of every critical point.
A critical point is called \emph{nondegenerate}
if its Hessian is nondegenerate.
A function with nondegenerate critical
points only is a \emph{Morse function}.
If $\Ss_V$ is Morse,
the change of topology of the sublevel set
$$
     \{\Ss_V\le a\}:=\{x\in\Ll M\mid
     \Ss_V(x)\le a\}
$$
when $a\in\R$ passes 
through a critical value 
is the subject of classical
Morse theory leading
to a CW-complex homotopy 
equivalent to $\{\Ss_V\le a\}$
(see e.g. Milnor~\cite{M64}).
In the case $V=0$
we use the notation
$\Ll^a M:=\{\Ss_0\le a\}$.

\subsubsection*{\boldmath$W^{1,2}$ Morse homology}
\label{subsubsec:W^{1,2}-HM}

A geometric reincarnation of
the idea of encoding the topology
of a sublevel set in terms of a Morse function
came (back) to light in 1982 through the
work of Witten~\cite{Wi82}.
Roughly speaking, the
\emph{Morse-Witten complex}
consists of chain groups
generated by the critical points of
a Morse function and a boundary
operator which counts flow lines
of the negative gradient flow
between critical points of
Morse index difference one
(for details see
e.g.~\cite{Sch93} and~\cite{JOA5}).
In recent years
Abbondandolo and Majer~\cite{AM04}
extended the theory from finite dimensions
to a class of Hilbert manifolds.
The free loop space fits into this framework
after completion with respect to the
Sobolev $W^{1,2}$ norm.
From now on we assume that $\Ss_V $
is a Morse function.
(A proof that this holds for a generic
potential $V$ is given in~\cite{JOA2}).
Then the set
$$
     \Pp^a(V)
     :=\{x\in\Pp(V)\mid \Ss_V(x)\le a\}
$$
is finite for every real number $a$.
Also from now on we assume that $a$
is a regular value of $\Ss_V$.
The chain groups
are the free abelian groups
generated by $\Pp^a(V)$
and graded by the
Morse index, namely
\begin{equation}\label{eq:chain-group}
     C^a_k(V)
     =\bigoplus_{\stackrel
     {\scriptstyle x\in\Pp^a(V)}
     {\IND_V(x)=k}}\Z x,\qquad
     k\in\Z.
\end{equation}
Our convention is that the direct
sum over an empty set equals $\{0\}$.
The negative of the $W^{1,2}$ gradient
vector field induces a flow
on the ($W^{1,2}$ completion of the)
loop space whose unstable manifolds
are of finite dimension and
whose stable manifolds are of finite
codimension. Let us choose an orientation
of the unstable manifold
of every critical point.
If the Morse-Smale condition holds
(this means that all stable and
unstable manifolds intersect 
transversally),
then there are only
finitely many flow lines
between critical points of index difference
one. These are called
\emph{isolated flow lines}.
Each one inherits an orientation,
because it is the intersection
of an oriented and a cooriented
submanifold. 
Let the characteristic sign
of an isolated flow line
be $+1$ if the inherited
orientation coincides with
the one provided by the flow
and $-1$ else.
Counting isolated flow lines with 
characteristic signs
defines the Morse-Witten boundary operator.
The associated homology groups
$\HM_*^a(\Ll M,\Ss_V,W^{1,2})$
are called \emph{$W^{1,2}$ Morse homology}.
By the theory of Abbondandolo and Majer
it is naturally isomorphic to
integral singular homology
$\Ho_*(\{\Ss_V\le a\})$.
If $\Ll M$ has several connected
components $\Ll_\alpha M$,
then there is a separate isomorphism
for each of them.
The label $\alpha$ denotes
a homotopy class
of free loops in $M$.

\subsubsection*{\boldmath$L^2$ Morse homology}
\label{subsubsec:L^2-HM}

Replacing the $W^{1,2}$ metric on
the free loop space by the $L^2$ metric
leads to a new boundary
operator on the chain
groups~(\ref{eq:chain-group}), namely
by counting negative 
$L^2$ gradient 'flow lines'.
In fact the $L^2$ metric
gives rise only to a semiflow
in forward time and so we
view -- in the
spirit of Floer theory -- 
the negative gradient 
flow equation on the loop
space as a PDE for 
smooth cylinders in $M$.
Flow lines are then replaced by
solutions $u:\R\times S^1\to M$
of the \emph{heat equation}
\begin{equation}\label{eq:heat}
     \p_s u-\Nabla{t}\p_t u-\nabla V_t(u)=0
\end{equation}
which satisfy
\begin{equation}\label{eq:heat-lim}
     \lim_{s\to\pm\infty} u(s,t) 
     = x^\pm(t),\qquad
     \lim_{s\to\pm\infty}\p_su(s,t)
     = 0.
\end{equation}
Here $x^\pm\in\Pp(V)$ and
the limits are supposed to be
uniform in $t$. 
Stable manifolds can still
be defined via the forward semiflow,
whereas to define
unstable manifolds we
use the heat flow lines~(\ref{eq:heat}).
The former are of
finite codimension
and the latter of finite
dimension. Hence
characteristic signs
can be assigned to isolated
(index difference one)
flow lines just as
in the case of $W^{1,2}$
Morse homology above.
The parabolic moduli space
$\Mm^0(x^-,x^+;V)$ is the
set of solutions of~(\ref{eq:heat})
and~(\ref{eq:heat-lim}).
In this setting we say that the
\emph{Morse-Smale condition} holds,
if the linear operator obtained by
linearizing~(\ref{eq:heat}) at a
solution $u$ is onto
for all $u\in \Mm^0(x^-,x^+;V)$
and all $x^\pm\in\Pp(V)$.
In this case $\Mm^0(x^-,x^+;V)$
is a smooth manifold whose dimension
equals the difference of the Morse indices.
In the case of index difference one
the quotient by the free time shift
action is a finite set.
Counting its elements
with characteristic signs
defines the $L^2$
Morse-Witten boundary
operator $\p^M$.
The associated homology 
$\HM_*^a(\Ll M,\Ss_V,L^2)$ is
called $L^2$ Morse homology.
It is naturally isomorphic
to $\Ho_*(\{\Ss_V\le a\})$.
(We should emphasize that this
is work in progress~\cite{JOA-FUTURE}).

\subsubsection*{Floer homology}
\label{subsubsec:HF}

The critical points of $\Ss_V$
can be interpreted  via the
Legendre transformation as the 
critical points of the
\emph{symplectic action functional}
$$
     \Aa_V(z)=\Aa_{H_V}(z)
     := \int_0^1 \biggl(
       \inner{y(t)}{\p_t x(t)} 
       - H_V(t,x(t),y(t))\biggr)\,dt.
$$
Here $z=(x,y)$ where
$x:S^1\to M$ is a smooth map
and $y(t)\in T_{x(t)}^*M$
depends smoothly on $t\in S^1$.
The Hamiltonian 
$H_V:S^1\times T^*M\to\R$ is
given by
$$
     H_V(t,q,p) = \frac12|p|^2+V_t(q).
$$
A loop
$z=(x,y)$ in $T^*M$ 
is a critical point of $\Aa_V$ iff $x$
is a critical point of $\Ss_V$
and $y(t)\in T_{x(t)}^*M$ is related to
$\p_t x (t)\in T_{x(t)}M$ via the isomorphism
$g:TM\to T^*M$ induced by the Riemannian metric. 
For such loops $z$ the symplectic action $\Aa_V(z)$ 
agrees with the classical action $\Ss_V(x)$.

In contrast to the classical
action, the symplectic action
is in general neither bounded
below nor do the critical
points admit finite Morse indices,
and most importantly
its $L^2$ gradient does not
define a flow on the loop space.
It was a great achievement
of Floer~\cite{F89} to
nevertheless set up
a Morse-Witten type complex.
His key idea was to reinterpret
the negative $L^2$ gradient equation
as elliptic PDE for maps from the cylinder
to the symplectic manifold imposing
appropriate boundary conditions
to make the problem Fredholm.
Floer's original setup
was a \emph{closed} symplectic manifold
subject to two topological assumptions
to ensure compactness of moduli spaces
and existence of a natural grading.
A Hamiltonian function $H$
and an almost complex structure $J$
need to be chosen 
to define the Floer complex.
The power of Floer theory
lies in the fact that Floer homology
is independent of these choices.
This is called the 
\emph{Floer continuation principle}.
Floer showed that if
$H$ is autonomous
and a $C^2$-small
Morse function,
then the Floer chain complex
equals the Morse-Witten complex.
Hence Floer homology is naturally
isomorphic to singular integral homology
of the closed symplectic manifold itself.
For introductory reading
we refer to Salamon's lecture notes~\cite{S97}
and the recent survey by Laudenbach~\cite{L04}.
A discussion
on a more advanced level, also
including applications
of Floer theory, can be found
in Chapter~12 of~\cite{MS04}.

Now consider the
cotangent bundle $T^*M$ equipped
with its canonical symplectic structure
$\omega_0=-d\theta$, where $\theta$
is the Liouville form. Since the
symplectic form is exact and
the first Chern class of $T^*M$
with respect to the  metric
induced almost complex structure
vanishes, both topological
assumptions of Floer are met. 
The former excludes
existence of nonconstant
$J$-holomorphic spheres,
which is an obstruction towards
compactness of the moduli spaces,
and the latter implies that
the Conley-Zehnder index
of 1-periodic Hamiltonian
orbits is well defined.
Since the 1-periodic orbits
of the Hamiltonian flow
are precisely the critical
points of $\Aa_V$ and these
coincide with the critical points
of $\Ss_V$ up to 
natural identification,
the Floer chain groups are again
given by the free abelian
groups~(\ref{eq:chain-group}).
(The standard Floer grading
is the negative Conley-Zehnder index,
which is proved in~\cite{JOA2} to
equal the Morse index; up to
a constant if $M$ is not orientable).
The Riemannian metric on $M$
provides the isomorphism
\begin{equation}\label{eq:splitting}
     T_{(x,y)}T^*M\to T_xM\oplus T_x^*M,
\end{equation}
which takes the derivative of a curve 
$\R\to T^*M:t\mapsto z(t)=(x(t),y(t))$
to the derivatives
of the two components, namely
$$
     \p_t z(t)
     \mapsto
     (\p_t x(t),\Nabla{t}y(t)).
$$
The metric also induces
an almost complex
structure $J_g$ and a metric $G_g$ on $T^*M$.
These and the symplectic
form are represented by
$$
     J_g=\begin{pmatrix}
     0&-g^{-1}\\g&0\end{pmatrix},\qquad
     G_g=\begin{pmatrix}
     g&0\\0&g^{-1}\end{pmatrix},\qquad
     \omega_0=\begin{pmatrix}
     0&-\1\\\1&0\end{pmatrix}.
$$
These three structures
are compatible in the sense that
$\omega_0(\cdot,J_g\cdot)=G_g(\cdot,\cdot)$.
Flow lines are then 
replaced by solutions
$w:\R\times S^1\to T^*M$
of Floer's equation
$$
     \p_s w+J_g(w)\p_t w-\nabla H_V(t,w)=0.
$$
It is the negative $L^2$ gradient
equation for the symplectic action
viewed as an elliptic PDE.
Its solutions are called
\emph{Floer trajectories} or
\emph{Floer cylinders}.
If we identify $T^*M$ and $TM$
via the metric isomorphism
and represent Floer's equation
with respect to the
splitting~(\ref{eq:splitting}),
we obtain the
pair of equations
\begin{equation}\label{eq:floer}
     \p_su-\Nabla{t}v-\nabla V_t(u)
     = 0,\qquad
     \Nabla{s}v+\p_tu-v=0
\end{equation}
for $(u,v):\R\times S^1\to TM$.
The Floer moduli space
$\Mm^1(x^-,x^+;V)$ is the set
of solutions
of~(\ref{eq:floer}) subject
to the boundary conditions
\begin{equation}\label{eq:floer-lim}
     \lim_{s\to\pm\infty}u(s,t)
     =x^\pm(t),\qquad
     \lim_{s\to\pm\infty}v(s,t)
     =\p_t x^\pm(t),
\end{equation}
and $\p_su$ and $\Nabla{s}v$ 
converge to zero
as $|s|\to\infty$, and all limits
are uniformly in $t$.
If the Morse-Smale condition
is satisfied,
then $\Mm^1(x^-,x^+;V)$ is
a smooth manifold whose dimension
is given by the difference
of Morse indices.
Since $T^*M$ is
noncompact, we
do not obtain for free
uniform apriori
$C^0$-bounds for the
Floer solutions as
in the standard case
of a \emph{closed} symplectic manifold.
Such bounds were established in 1992
by Cieliebak in his diploma thesis
(published in~\cite{C94})
and recently extended to
a class of radial Hamiltonians
by the author~\cite{JOA4}
and to another class
of not necessarily radial
Hamiltonians by Abbondandolo
and Schwarz~\cite{AS04}.
Given these bounds,
proving compactness 
and setting up Floer homology
is standard.
In the case that $\Ll T^*M$ has several
connected components,
Floer homology is defined for
each of them separately and
denoted by
$\HF_*^a(T^*M,H_V,J_g;\alpha)$.
Here $\alpha$ denotes a homotopy
class of free loops in $M$.
(Throughout we identify
homotopy classes of free loops
in $M$ and in $T^*M$).
In fact Floer homology can be
defined for classes of
Hamiltonians more general
than $H_V$, for instance
those growing
quadratically in $p$
outside of a compact set
(see Section~\ref{sec:AS})
or convex radial Hamiltonians
(see Section~\ref{sec:noncontractible}).

\subsubsection*{Isomorphisms 
                between the theories}
\label{subsubsec:relations}

When Floer homology for
physical Hamiltonians of the
type kinetic plus potential energy 
could be defined by~1992
due to Cieliebak's breakthrough,
the obvious question was
\emph{`What is it equal to?'}.

The first
%
answer\footnote{By then
Theorem~\ref{thm:main}
had been expected to be true by
part of the community.
For instance the problem was proposed
as a PhD project to the present author
by Helmut Hofer during
winter term 1993/94 at ETH Z\"urich.
In summer~1996
upon meeting Dietmar Salamon in Oberwolfach
we matched up and started our joint approach.
A short time later~\cite{V96} appeared.
In private communication at a Warwick conference,
around~1998, Matthias Schwarz first told me about
an alternative approach via a
mixed boundary value problem.}
%
in the literature
is due to Viterbo who
conjectured in his~1994 
ICM talk~\cite{V95} -- based on
formal interpretation of the symplectic
action functional as a generating function --
that Floer homology of the cotangent bundle
represents singular homology of the
free loop space.
In his 1996 preprint~\cite{V96}
he gave a beautiful line of argument
(for the component $\Ll_0 M$
of contractible loops, case $V=0$,
coefficients in $\Z_2$).
The idea is to
view the time-1-map
$\varphi_1$ of the Hamiltonian flow
as an $r$-fold composition
of symplectomorphisms close
to the identity
(set $\psi:=\varphi_{1/r}$)
in order to arrive at the well known
finite dimensional approximation
of the free loop space
via broken geodesics.
While the argument
consists of numerous steps and the idea
of each one is described in detail,
not all technical details are
provided.

Also in 1996 a first example
was computed by the present
author~\cite{JOA0c}, namely
Floer homology of the cotangent
bundles of the euclidean torus
$\T^n=\R^n/2\pi\Z^n$
confirming the conjecture
for all connected components
of $\Ll\T^n$, for every $n\in\N$.
This is reviewed in
Section~\ref{sec:torus}.

In 2003 Salamon and the
present author~\cite{JOA3}
proved existence of a natural
isomorphism
$$
     \HM_*^a(\Ll M,\Ss_V,L^2;\Z)
     \to
     \HF_*^a(T^*M,H_V,J_g;\Z).
$$
(Partial results were established
in~1999 in the
PhD thesis~\cite{JOA1}).
The idea is to introduce a
real parameter $\eps>0$ and to replace
the standard almost complex structure $J_g$
by $J_{\eps^{-1}g}$.
Both Floer homologies are naturally isomorphic
by Floer continuation.
The key step is then
to prove that for every
sufficiently small $\eps$
the parabolic and elliptic moduli
spaces can be identified.
This means that
the $\eps$-Floer
and the $L^2$ Morse chain complexes
are \emph{identical}.

In 2004 Abbondandolo and 
Schwarz~\cite{AS04}
proved existence of a natural isomorphism
$$
      \HM_*^a(\Ll M,\Ss_V,W^{1,2})
      \to
      \HF_*^a(T^*M,H_V,J_g)
$$
by constructing
a \emph{chain isomorphism}
in the case of orientable $M$.
Their approach works
for more general Hamiltonians
and almost complex structures
(see Section~\ref{sec:AS}).
The idea is to study
a mixed boundary value problem for
Floer half cylinders
$w=(u,v):[0,\infty)\to T^*M$.
For $s\to+\infty$ the standard Floer
boundary condition~(\ref{eq:floer-lim})
is imposed, whereas at $s=0$
the base loop $u(0,\cdot)$
is required to belong to an unstable
manifold of the negative $W^{1,2}$
gradient flow of $\Ss_V$.

Also in~2004 the present
author~\cite{JOA4} extended
the definition and
computation of Floer homology 
to the class of convex
radial Hamiltonians
(those of the form
$h=h(\abs{p})$ with $h^{\prime\prime}\ge 0$).

The main result of~\cite{V96},
\cite{JOA3} and~\cite{AS04}
as formulated in~\cite{JOA3} is
the following theorem.

\begin{theorem}\label{thm:main}
Let $M$ be a closed Riemannian
manifold.
Assume $\Ss_V$ is Morse and $a$
is either a regular
value of $\Ss_V$ or is equal to infinity. 
Then there is a natural isomorphism
$$
     \HF^a_*(T^*M,H_V,J_g;R)
     \simeq\mathrm{H}_*(\{\Ss_\Vv\le a\};R)
$$
for every principal ideal domain $R$. 
If $M$ is not simply connected,
then there is a separate
isomorphism for each component of the loop space. 
The isomorphism commutes with
the homomorphisms
$\HF^a_*(T^*M,H_V,J_g)\to\HF^b_*(T^*M,H_V,J_g)$
and $\mathrm{H}_*(\{\Ss_\Vv\le a\})
\to\mathrm{H}_*(\{\Ss_\Vv\le b\})$, for $a<b$,
which are induced by inclusion.
\end{theorem}

We summarize the discussion
by the diagram below
in which arrows represent isomorphisms.
The homologies are defined as usual
by first perturbing to achieve Morse-Smale
transversality and then taking the
homology of the perturbed chain complex.
The branch on the right hand side
indicates Viterbo's
finite dimensional approximation
argument (which he actually formulated
in terms of cohomology;
see Section~\ref{sec:V}).
\begin{equation*}
\begin{split}
\xymatrix{
     \HF_*^a(T^*M,H_0,J_{\eps^{-1}g})
    &
     \HF_*^a(T^*M,H_0,J_g)
     \ar[r]_{\textstyle\cite{V96}}
     \ar[l]^{\overset{\scriptstyle\text{Floer}}
             {\scriptstyle\text{continuation}}}
    &
     \HF^a_*(\Delta_r,\Gamma_r(\varphi^{H_0}))
     \ar[d]^{\textstyle\cite{V96}}
     \\
     \HM_*^a(\Ll M,\Ss_0,L^2)
     \ar[dr]_{\textstyle\cite{JOA-FUTURE}}
     \ar[u]^{\textstyle\cite{JOA3}}
    &
     \HM_*^a(\Ll M,\Ss_0,W^{1,2})
     \ar[d]^{\textstyle\cite{AM04}}
     \ar[u]_{\textstyle\cite{AS04}}
    &
     \HM^a_*(U_{r,\eps},S_r)
     \ar[d]^{\textstyle\cite{V96,V97}}
     \\
    &
     \Ho_*(\Ll^a M)
    &
     \Ho_*(\Lambda^a_r)
     \ar[l]_{\textstyle\cite{V97}}
}
\end{split}
\end{equation*}
It would be interesting
to fill in the missing link between
$L^2$ and $W^{1,2}$ Morse
homology, i.e. construct an
isomorphism which is natural
in the sense that the
corresponding triangle
and rectangle in the diagram
are both commutative.

The remaining part of this
text is organized as follows.
We present the three appoaches
towards Theorem~\ref{thm:main}
in chronologically reverse order
in Sections~\ref{sec:AS}--\ref{sec:V}.
This way complexity increases 
-- as it should be.
In Section~\ref{sec:torus} we calculate
Floer homology of the cotangent
bundle of the euclidean $n$-torus.
An application of
Theorem~\ref{thm:main} to
existence of noncontractible
periodic orbits is reviewed in
Section~\ref{sec:noncontractible}.
For an application towards
Arnold's chord conjecture
we refer to Cieliebak's
paper~\cite{C02}.

\medskip\noindent
{\small\it Acknowledgements:}
{\footnotesize
We gratefully acknowledge
partial financial support
by DFG SPP~1154
{\sl Globale Differentialgeometrie},
Scuola Normale Superiore Pisa
and \'Ecole Polytechnique Paris.
We are particularly indebted to
Alberto Abbondandolo, Kai Cieliebak,
Pietro Majer and
Claude Viterbo
for numerous helpful and pleasant
conversations on the subject.}

\section{Mixed boundary value problem}
\label{sec:AS}

In three steps we review
the approach of Abbondandolo
and Schwarz~\cite{AS04}.
They assume for simplicity that
$M$ is orientable.

First of all, the authors
set up Floer theory for a more
general class of Hamiltonians $H$
and almost complex structures $J$.
Let $\Pp(H)$ denote the set of
critical points of the symplectic
action $\Aa_H$. These are precisely the
1-periodic orbits of the
Hamiltonian flow on $T^*M$.
The crucial (metric independent)
assumptions on $H$
are the following.
Outside of a compact set
$H$ is supposed to satisfy
$$
     dH(t,q,p) \; p\p_p-H(t,q,p)
     \ge c_0\abs{p}^2-c_1,
     \eqno{\text{(H1)}}
$$
$$
     \Abs{\Nabla{q} H(t,q,p)}
     \le c_2(1+\abs{p}^2),\qquad
     \Abs{\Nabla{p} H(t,q,p)}
     \le c_2(1+\abs{p}),
     \eqno{\text{(H2)}}
$$
for some constants $c_0>0$
and $c_1,c_2\ge 0$.
Assumptions~(H1) and~(H2)
guarantee that the set $\Pp^a(H)$
is finite for every real number $a$,
whenever $\Aa_H$ is Morse (which we shall
assume from now on, since
it is true for generic $H$).
More importantly, assumptions~(H1) and~(H2)
allow the authors to establish
$C^0$-bounds for Floer solutions
associated to
almost complex structures $J$
sufficiently $C^\infty$-close to $J_g$.
Then the definition of the chain complex
$(\CF_*(H),\p_*(H,J))$
proceeds by standard arguments.
Any two choices of $J$
lead to isomorphic \emph{chain complexes}.
In contrast homology is independent of $(H,J)$
and denoted by $\HF_*(T^*M)$.

Secondly, the Morse-Witten complex is defined
for the Hilbert manifold $W^{1,2}(S^1,M)$
and the classical action functional
$\Ss_L(x):=\int_0^1 L(t,x(t),\dot x(t))\:dt$.
Here the admissible Lagrangians $L$
are those for which
there exist constants
$d_0>0$ and $d_1\ge 0$ such that
$$
     \Nabla{vv} L(t,q,v)\ge d_0\1,
     \eqno{\text{(L1)}}
$$
\begin{minipage}{11cm}
$$
     \Abs{\Nabla{qq} L(t,q,v)}
     \le d_1(1+\abs{v}^2),\qquad
     \Abs{\Nabla{qv} L(t,q,v)}
     \le d_1(1+\abs{v}),
$$
$$
     \Abs{\Nabla{vv} L(t,q,v)}
     \le d_1,
$$
\end{minipage} 
\hfill
\begin{minipage}{.66cm}
(L2)
\end{minipage}

\vspace{.2cm}
\noindent
for all $(t,q,v)\in S^1\times TM$.
Perturbing $L$
if necessary we assume
from now on that $\Ss_L$
is Morse and denote
the set of its critical points
by $\Pp(L)$.
The classical action
exhibits a number of
rather useful features.
For instance, it satisfies
the Palais-Smale condition,
it is bounded from below
and its critical points have
finite Morse indices (which
equal minus the corresponding
Conley-Zehnder indices)\footnote{The
sign difference
between~\cite{JOA2,JOA3} and~\cite{AS04}
is due to the different normalizations
$\mu_{CZ}(t\mapsto e^{i\pi t})=1$
and $\mu_{CZ}(t\mapsto e^{-i\pi t})=1$
(with $t$ running through $[0,1]$),
respectively.}.
Choosing an auxiliary
Morse-Smale metric
$\Gg$ on the Hilbert manifold,
the work of Abbondandolo and
Majer~\cite{AM04} establishes
existence of the Morse complex
$(\CM_*(\Ss_L),\p_*(\Ss_L,\Gg))$
and shows that its homology
is naturally isomorphic to the singular
homology of the free loop space.

Given $L$ and $\Gg$ as in Step~2,
the crucial third step 
is to construct
a grading preserving chain complex
isomorphism
$$
     \Theta_*:
     (\CM_*(\Ss_L),\p_*(\Ss_L,\Gg))
     \to
     (\CF_*(H_L),\p_*(H_L,J))
$$
where the Hamiltonian $H_L$
arises from the Lagrangian $L$
via Legendre transformation
$(t,q,p)\mapsto (t,q,v(t,q,p))$.
More precisely, define
$$
     H_L(t,q,p)
     :=\max_{v\in T_qM}\left(
     \langle p,v \rangle-L(t,q,v)\right).
$$
For each $(t,q,p)$ the maximum
is achieved at a unique point
$v(t,q,p)$ by condition~(L1).
The Legendre transformation
provides a natural identification
of the critical points
of $\Aa_{H_L}$ and $\Ss_L$, namely
$\Pp(H_L)\to\Pp(L):
(x,y)\mapsto x$.
This shows that both chain groups coincide.
\\
To define a chain homomorphism
fix $q\in\Pp(L)$ and
$z\in\Pp(H_L)$. Then
consider half cylinders
$w:[0,\infty)\times S^1\to T^*M$
solving Floer's equation
$$
     {\overline\p}_{J,L} w
     :=\p_sw+J(w)\p_t w
     -\nabla H_L(t,w)=0
$$
and such that $w(s,\cdot)$
converges uniformly to $z$,
for $s\to+\infty$.
The boundary condition
at the other end is
that the loop
$w(0,\cdot)$ projects to the 
unstable manifold of $q$, i.e.
$u:=\pi\circ w(0,\cdot) \in W^u(q)$
where $\pi:T^*M\to M$ is
the bundle projection.
These half cylinders are the
elements of the moduli space
\begin{equation*}
\begin{split}
     \Mm^+(q,z)
    :=
     \{ w:[0,+\infty)\times S^1\to T^*M
     \mid\:
    &\text{$\pi\circ w(0,\cdot)\in W^u(q)$,
     ${\overline\p}_{J,L} w=0$,}\\
    &\text{$\lim_{s\to\infty} 
     w(s,\cdot)=z(\cdot)$}
     \}.
\end{split}
\end{equation*}
The problem is Fredholm,
because the unstable manifolds $W^u(q)$
are finite dimensional
and the boundary conditions for
the $\overline\p$-equation are Lagrangian
and nondegenerate, respectively.
For generic $J$ this
moduli space is a smooth manifold
of dimension $\IND_L(q)-\IND_L(z)$.
In the case of equal indices
it is a discrete set and $\Theta_*$
is defined by counting its elements.
Of course, compactness of
the moduli space needs to be
established first.
Here a crucial observation of Abbondandolo
and Schwarz enters, namely the inequality
\begin{equation}\label{eq:AS-inequality}
     \Aa_{H_L}(x,y)\le\Ss_L(x)
\end{equation}
for every loop $(x,y):S^1\to T^*M$.
On critical points equality holds.
To prove compactness
a uniform action bound is needed
to start with, but this follows
immediately from~(\ref{eq:AS-inequality}):
if $w\in\Mm^+(q,z)$, then
$\Aa_{H_L}(w(s,\cdot))\le\Ss_L(q)$
for all $s\ge0$.
Compactness now follows from
the $C^0$ estimate discussed
in the first step.
A gluing argument proves that
$\Theta_*$ is a chain map.
Moreover,
inequality~(\ref{eq:AS-inequality})
shows that
$\Mm^+(q,z)=\emptyset$
whenever $\Ss_L(q)\le\Aa_{H_L}(z)$,
unless $q$ and $z$ correspond
to the same critical point.
In this case $\Mm^+(q,z)$
consists of a single element,
the constant solution
$w(s,\cdot)=z(\cdot)$.
There is a differential version
of~(\ref{eq:AS-inequality}) 
at $z\in\Pp(H_L)$, namely
$$
     d^2\Aa_{H_L}(z)\:(\cdot,\cdot)
    \le d^2\Ss_L(\pi(z))\:
     (D\pi(z)\cdot,D\pi(z)\cdot).
$$
It is used to prove that
at a constant
solution $w(s,\cdot)=z(\cdot)$
the Morse-Smale
condition is automatically true.
Ordering the generators
of the chain groups
by increasing action,
it follows that $\Theta_k$
is an upper triangular matrix
with diagonal entries $\pm1$.
This proves that
$\Theta_*$ is a chain
isomorphism.

\section{Singular perturbation
         and adiabatic limit}
\label{sec:SW}

This section reviews
the approach by Salamon and the
present author who
established a natural isomorphism
$\HF_*^a(T^*M,H_V,J_g)\to
\HM_*^a(\Ll M,\Ss_V,L^2)$
in~\cite{JOA3}.
It is work in
progress~\cite{JOA-FUTURE}
to show that the latter
is naturally isomorphic to
$\Ho_*(\Ll^a M)$.
The main idea to relate Floer and
Morse homology is to
introduce a real parameter $\eps>0$
in Floer's equations:
replace $(J_g,G_g)$ by 
$$
     (J_\eps,G_\eps)
     :=(J_{\eps^{-1}g},G_{\eps^{-1}g}).
$$
By Floer's continuation principle
this will not change Floer homology.
Let us identify $T^*M$ with $TM$
via the metric isomorphism. A map
$(\tilde u,\tilde v):\R\times S^1\to TM$
is a solution to the
$(J_\eps,G_\eps)$-Floer
equations~(\ref{eq:floer})
if and only if
$u(s,t):=\tilde u(\eps^{-1}s,t)$
and $v(s,t):=v(\eps^{-1}s,t)$
satisfy
\begin{equation}\label{eq:floer-eps}
     \p_su-\Nabla{t}v-\nabla V_t(u)
     = 0,\qquad\Nabla{s}  v
     +\eps^{-2}\left(\p_tu-v\right)=0.
\end{equation}
Denote the space of solutions
to~(\ref{eq:floer-eps})
and~(\ref{eq:floer-lim}) by
$\Mm^\eps(x^-,x^+;V)$.
We shall outline how to
prove that there is
a one-to-one correspondence
between the solutions
of~(\ref{eq:floer-eps})
and~(\ref{eq:heat})
subject to boundary
conditions~(\ref{eq:floer-lim})
and~(\ref{eq:heat-lim}), respectively,
whenever the index difference
is one\footnote{The case of
arbitrary index difference
is closely related to Cohen's
conjecture~\cite{Co04}. It
is the missing link
in proving that cylindrical
Gromov-Witten invariants
of the cotangent bundle
represent string topology
of the free loop space.
A different proof relating
the particular cases of
the three point invariant
provided by the pair of pants product
in Floer homology
and the Chas-Sullivan
loop product~\cite{CS99} 
is in preparation
by Abbondandolo and
Schwarz~\cite{AS04b}.
}.
This then shows that both
chain complexes are identical.
A first hint that such a
bijection between parabolic and
$\eps$-elliptic flow lines
exists is provided
by the \emph{energy identity}
\begin{equation*}
\begin{split}
     E^\eps(u,v)
    &:=
     \frac12\int_{-\infty}^\infty
     \int_0^1\left(
     |\p_su|^2
     +|\Nabla{t}v+\nabla V_t(u)|^2
     +\eps^2|\Nabla{s}v|^2
     +\eps^{-2}|\p_tu-v|^2
     \right) \\
    &=
     \Ss_V(x^-)-\Ss_V(x^+)
\end{split}
\end{equation*}
for the solutions
of~(\ref{eq:floer-eps})
and~(\ref{eq:floer-lim}).
It shows that $\p_tu-v$ must
converge to zero in
$L^2$ as $\eps\to0$, but if
$\p_tu=v$ then the first equation
in~(\ref{eq:floer-eps}) is equivalent
to~(\ref{eq:heat}).

The next step would be to prove
that $\Ss_V$ is Morse-Smale
for generic $V$.
Unfortunately, this remains
an \emph{open problem}.
Instead we introduce a more
general class of perturbations
for which Morse-Smale transversality
can be achieved generically
by standard methods.
These perturbations take
the form of smooth functions
$\Vv:\Ll M\to\R$ satisfying
a list of axioms (which
contains the properties used at some
point in the proof; see~\cite{JOA3}).
Assume for the moment that $M$ was
embedded isometrically in some euclidean
space $\R^N$, fix a loop $x_0$
in $M$ and let $\rho:\R\to[0,1]$
be a smooth cutoff function.
Then a typical example of an
abstract perturbation is given by
$$
     \Vv(x)
     :=\rho\left(\left\|x
     -x_0\right\|_{L^2}^2\right)
     \int_0^1V_t(x(t))\,dt.
$$
With new functionals
$\Ss_\Vv:=\Ss_0+\Vv$
and $\Aa_\Vv:=\Aa_0+\Vv$
equation~(\ref{eq:floer-eps})
turns into
\begin{equation}\label{eq:floer-V}
     \p_su-\Nabla{t}v-\grad\Vv(u)
     =0,\qquad
     \Nabla{s}v+\eps^{-2}(\p_tu-v)
     =0
\end{equation}
and the limit equation is
of the form
\begin{equation}\label{eq:heat-V}
     \p_su-\Nabla{t}\p_tu-\grad\Vv(u)
     =0.
\end{equation}
Here the $L^2$ gradient
$\grad\Vv(x)\in\Om^0(S^1,x^*TM)$
of $\Vv$ at $x\in\Ll M$ is defined by 
$$
     \int_0^1\inner{\grad\Vv(u)}
     {\p_su}\,dt 
     := \frac{d}{ds}\Vv(u)
$$
for every smooth path
$\R\to\Ll M:s\mapsto u(s,\cdot)$. 
The set
$\Pp(\Vv)$ consists of loops
$x:S^1\to M$ satisfying
$\Nabla{t}\p_t x=-\grad\Vv(x)$.
Define $\Mm^\eps(x^-,x^+;\Vv)$ and
$\Mm^0(x^-,x^+;\Vv)$
as before with~(\ref{eq:floer-V})
and~(\ref{eq:heat-V})
replacing~(\ref{eq:floer-eps})
and~(\ref{eq:heat}),
respectively.

Let $V_t$ be a potential such that
$\Ss_V$ is a Morse function and denote
$$
\Vv(x):=\int_0^1V_t(x(t))\,dt.
$$ 
Observe that this choice reproduces
$\Ss_V$ and $\Aa_V$, hence the
geometric equations~(\ref{eq:heat})
and~(\ref{eq:floer}).
Fix a regular value $a$ of $\Ss_V$ and
choose a sequence of perturbations
${\Vv_i:\Ll M\to\R}$
converging to $\Vv$ in the 
$\Cinf$ topology and such that
$\Ss_{\Vv_i}:\Ll M\to\R$
is Morse--Smale for every $i$.
We may assume without 
loss of generality that the
perturbations agree with 
$\Vv$ near the critical points
and that $\Pp(\Vv_i)=\Pp(V)$
for all $i$. Assume there
is a sequence $\eps_i>0$
converging to zero such that,
for every $\eps_i$ and
every pair $x^\pm\in\Pp^a(V)$
of index difference one,
there is a ($s$-shift equivariant)
bijection
\begin{equation}\label{eq:bijection}
     \Tt^{\eps_i}:
     \Mm^0(x^-,x^+;\Vv_i)\to
     \Mm^{\eps_i}(x^-,x^+;\Vv_i).
\end{equation}
The following diagram --
in which arrows represent isomorphisms --
shows how this implies our goal
(with $\Z_2$-coefficients).
\begin{equation*}
\begin{split}
\xymatrix{
     \HF_*^a(T^*M,\Vv_i,J_{\eps_i})
     \ar[d]
          ^{(\ref{eq:bijection})}
    &
     \HF_*^a(T^*M,H_V+W,J_{\eps_i})
     \ar[l]_{\text{Floer}\quad}
           ^{\text{contin.}\quad}
    &
     \HF_*^a(T^*M,H_V,J_{\eps_i})
     \ar[l]_{\quad\;\text{def}}
     \\
     \HM_*^a(\Ll M,\Ss_{\Vv_i},L^2)
     \ar[d]^{\textstyle\cite{JOA-FUTURE}}
    &
    &
     \\
     \Ho_*(\{\Ss_{\Vv_i}\le a\})
     \ar[rr]^{\text{homotopy
                    equivalent spaces}}
    &
    &
     \Ho_*(\{\Ss_V\le a\})
}
\end{split}
\end{equation*}
Starting at the upper right corner,
the first step is by
definition of Floer homology
for nonregular $H_V$, namely
add a small Hamiltonian
perturbation $W$ making
$H_V+W$ regular.
Floer continuation shows
independence of the choice.
Again by Floer's continuation argument
small Hamiltonians
and small abstract perturbations
lead to isomorphic homology groups.

It remains to construct the bijection
$\Tt^\eps$ in the case of index
difference one.
Assume throughout that
$\Ss_\Vv$ is Morse-Smale and fix a regular
value $a$ and a pair $x^\pm\in\Pp^a(V)$
of index difference one.
Denote $\Mm^0:=\Mm^0(x^-,x^+;\Vv)$ and 
$\Mm^\eps:=\Mm^\eps(x^-,x^+;\Vv)$.

\subsubsection*{Existence and uniqueness}
\label{subsubsec:exist-unique}

Given a parabolic solution we shall prove
by Picard-Newton iteration
existence and uniqueness of an elliptic
solution nearby.
The first step is to define, for smooth maps
$(u,v):\R\times S^1\to TM$ and
$p>2$, a map $\Ff^\eps_{u,v}$
between the Banach spaces
of $W^{1,p}$ and $L^p$ sections of the bundle
$u^*TM\oplus u^*TM\to\R\times S^1$.
To obtain uniform estimates for
small $\eps$ we introduce
weighted norms. The weights for $p=2$
are suggested by the energy identity.
Let $\Ff_\eps(u,v)$
be given by the left hand side of
the $\eps$-equations~(\ref{eq:floer-eps})
and denote by
$\Phi(x,\xi):T_xM\to T_{\exp_x(\xi)}M$
parallel transport along the
geodesic $\tau\mapsto\exp_x(\tau\xi)$. 
For compactly supported vector fields
$\zeta=(\xi,\eta)\in
\Om^0(\R\times S^1,u^*TM\oplus u^*TM)$
define the weighted norms
$$
     \left\|\zeta\right\|_{0,p,\eps}
     := \left(\int_{-\infty}^\infty\int_0^1
       \left(\left|\xi\right|^p 
       + \eps^p\left|\eta\right|^p\right)\,dtds
       \right)^{1/p},
$$
\begin{equation*}
\begin{split}
     \left\|\zeta\right\|_{1,p,\eps}
    &:=
     \biggl(\int_{-\infty}^\infty\int_0^1
       \bigl(\left|\xi\right|^p 
       + \eps^p\left|\eta\right|^p
       + \eps^p\left|\Nabla{t}\xi\right|^p 
       + \eps^{2p}\left|\Nabla{t}\eta\right|^p \\
    &\quad
       +\, \eps^{2p}\left|\Nabla{s}\xi\right|^p 
       + \eps^{3p}\left|\Nabla{s}\eta\right|^p
       \bigr)\,dtds\biggr)^{1/p},
\end{split}
\end{equation*}
and $\Ff_{u,v}^\eps:
W^{1,p}(\R\times S^1,u^*TM\oplus u^*TM) \to
L^p(\R\times S^1,u^*TM\oplus u^*TM)$ by
\begin{equation*}
     \Ff_{u,v}^\eps \begin{pmatrix} \xi \\
     \eta \end{pmatrix}
     :=
     \begin{pmatrix} \Phi(u,\xi)^{-1} & 0 \\
     0 & \Phi(u,\xi)^{-1} \end{pmatrix}
     \Ff_\eps \begin{pmatrix} \exp_u \xi \\
     \Phi(u,\xi) (v + \eta) \end{pmatrix}.
\end{equation*}
Abbreviate $\Dd_{u,v}^\eps
:=d\Ff_{u,v}^\eps(0,0)$ and
$\Dd_u^\eps:=\Dd_{u,\p_tu}^\eps$.

Let us fix a parabolic cylinder
$u\in\Mm^0$,
viewed as an approximate
solution of the $\eps$-elliptic
equations. Equivalently,
we view the origin as
approximate zero
of the map $\Ff_u^\eps
:=\Ff_{u,\p_tu}^\eps$
between Banach spaces.
Carrying out the Newton-Picard iteration
for this map we shall prove existence
of a true zero nearby.
The iteration method works if,
firstly, the initial value is small,
secondly, the linearized operator admits
a right inverse, and thirdly,
second derivatives
of the map can be controlled.
These conditions must be satisfied
uniformly for small $\eps>0$.
Choosing the origin as the
initial point of the iteration
we observe that
$$
     \Norm{\Ff_u^\eps(0,0)}_{0,p,\eps}
     =\Norm{\Ff_\eps(u,\p_tu)}_{0,p,\eps}
     =\Norm{(0,\Nabla{s}\p_tu)}_{0,p,\eps}
     =\eps \Norm{\Nabla{s}\p_tu}_p
     \le c_0\eps.
$$
The second identity uses
the heat equation~(\ref{eq:heat})
and the final estimate is by
exponential decay of heat flow solutions
with nondegenerate boundary conditions.
Verification of the second condition relies
heavily on the fact that $\Dd_u^\eps$
is Fredholm and surjective. (Fredholm
follows by nondegeneracy of the boundary
conditions~(\ref{eq:floer-lim})
and surjectivity is a consequence
of the Morse-Smale assumption
for $\Ss_\Vv$).
Let ${\Dd_u^\eps}^*$
be the adjoint operator of $\Dd_u^\eps$
with respect to the $L^2$ inner product
$\langle\cdot,\cdot\rangle_\eps$ with
associated norm $\norm{\cdot}_{0,2,\eps}$.
A right inverse of $\Dd_u^\eps$
is given by
$
     \Qq_u^\eps
     :={\Dd_u^\eps}^*
     \left(\Dd_u^\eps{\Dd_u^\eps}^*\right)^{-1}.
$
It allows to solve the equation
$0=\Dd_u^\eps\zeta_0+\Ff_u^\eps(0,0)$
and provides the correction term
$\zeta_0:=-\Qq_u^\eps\Ff_u^\eps(0,0)$.
Recursively, for $\nu\in\N$, define
the sequence of correction terms 
$\zeta_\nu=(\xi_\nu,\eta_\nu)$
by
\begin{equation*}
     \zeta_\nu
    :=- {\Dd_u^\eps}^* 
     (\Dd_u^\eps{\Dd_u^\eps}^*)^{-1}
     \Ff_u^\eps (Z_\nu),\qquad
     Z_\nu
    :=\sum_{\ell=0}^{\nu-1} \zeta_\ell.
\end{equation*}
To estimate the terms $\zeta_\nu$
it is crucial to have an estimate
for $\Dd_u^\eps$ on the image of its adjoint
operator. This is a consequence
of a Calderon-Zygmund estimate
for a Cauchy-Riemann type operator.
It follows that $Z_\nu$
is a Cauchy sequence.
The third condition is verified
by quadratic estimates, i.e.
estimates for
$$
     \Ff_u^\eps(Z_\nu+\zeta_\nu)
     -\Ff_u^\eps(Z_\nu)
     -d\Ff_u^\eps(Z_\nu)\zeta_\nu
     ,\qquad
     d\Ff_u^\eps(Z_\nu)\zeta_\nu
     -\Dd_u^\eps\zeta_\nu.
$$
One proves by induction that
there is a constant $c>0$ such that
$$
     \left\|\Ff_u^\eps(Z_{\nu+1})
     \right\|_{0,p,\eps^{3/2}}
     \le\frac{c}{2^\nu}\eps^{7/2-3/2p}.
$$
for all $\nu\in\N$ and $\eps>0$ small.
Therefore the Cauchy sequence $Z_\nu$
converges to a zero $Z$ of $\Ff_u^\eps$
with $Z\in\im {\Dd_u^\eps}^*$
and $\Norm{Z}_{1,p,\eps}\le c\eps^2$.
Now $Z=(X,Y)$ corresponds to
a zero $(u^\eps,v^\eps)$ of $\Ff_\eps$,
more precisely to an element of $\Mm^\eps$,
and we define
$$
     \Tt^\eps(u)
     :=(u^\eps,v^\eps)
     :=(\exp_u X,\Phi(u,X)(\p_t u+Y)).
$$
Using again the quadratic
estimates, roughly speaking,
one can prove uniqueness
of the constructed solution
$(u^\eps,v^\eps)$ in an even
bigger neighbourhood of $(u,\p_tu)$.
Here a crucial assumption is
that the 'difference' $Z$ between
the parabolic and $\eps$-elliptic solution
is in the image of ${\Dd_u^\eps}^*$.
It follows that $\Tt^\eps$ is well defined.
Moreover, it is time shift equivariant.
Injectivity of $\Tt^\eps$ follows,
because the quotient
$\Mm^0/\R$ by the time shift action
is a finite set and so admits
a positive smallest distance
between its elements.
Since the existence and uniqueness
range shrinks like $\eps^2$, injectivity holds
for all sufficiently small $\eps$.

\subsubsection*{Surjectivity}
\label{subsubsec:surjectivity}

In four steps we sketch the proof that
the map $\Tt^\eps:\Mm^0\to\Mm^\eps$
is surjective whenever
$\eps>0$ is sufficiently
small.

The first step establishes
uniform \emph{apriori $L^\infty$ bounds}
for elliptic solutions
$(u^\eps,v^\eps)\in\Mm^\eps$ and
their first 
and second derivatives.
In each case this is based
on first proving
\emph{slicewise} $L^2$ bounds,
i.e. bounds in $L^2(S^1)$
for arbitrary fixed $s\in\R$,
and then applying a mean value inequality
for the operator $L_\eps:=\eps^2\p_s^2
+\p_t^2-\p_s$. More precisely, we prove that
there is a constant $c>0$ such that
for all $\eps>0$, $r\in(0,1]$, and
$\mu\ge0$ the following is true.
If $w$ is a smooth function
(for instance $w=\abs{v^\eps}^2$)
defined on the
\emph{parabolic domain} $P_r^\eps:=
(-r^2-\eps r,\eps r)\times(-r,r)$, then
$$
     L_\eps w
     :=\left(\eps^2{\p_s}^2+{\p_t}^2
       -\p_s\right)w
     \ge-\mu w, \quad w\ge0 
$$
imply
$$
     w(0)
     \le\frac{2ce^{\mu r^2}}{r^3}
     \int_{P_r^\eps}w.
$$
Now use the slicewise $L^2$
estimates on the right hand side
of the inequality to obtain
the $L^\infty$ estimate.
A bubbling argument
enters the proof in the case of
first derivatives,
as is expected by comparison
with standard Floer theory.
It provides a rather weak
(in terms of powers of $\eps$)
preliminary $L^\infty$
estimate which, however, suffices
to establish the slicewise $L^2$ bounds.

The second step proves \emph{uniform
exponential decay} of $\Abs{\p_su^\eps(s,t)}
+\Abs{\Nabla{s}v^\eps(s,t)}$
towards the ends of the cylinder.
The standard method of proof
works uniformly in $\eps\in(0,1]$
under the assumption that
the energy near the ends of the cylinder
is uniformly bounded by a small
constant (no energy concentration
near infinity).

In Step three we establish
\emph{local surjectivity} by means of
a time shift argument.
Given $u\in\Mm^0$
and $(u^\eps,v^\eps)\in\Mm^\eps$
sufficiently close to $(u,\p_tu)$,
define $\zeta=(\xi,\eta)$
by $u^\eps=\exp_u(\xi)$ and
$v^\eps=\Phi(u,\xi)(\p_tu+\eta)$.
The idea is to prove that
after a suitable time shift
the pair $\zeta=(\xi,\eta)$ 
satisfies the crucial hypothesis
$\zeta\in\im\,{\Dd^\eps_u}^*$
in the uniqueness theorem
and therefore
$(u^\eps,v^\eps)=\Tt^\eps(u)$.
More precisely, since the Fredholm
index of $\Dd^\eps_u$
is given by the Morse index
difference
and $\Dd^\eps_u$ is surjective
by the Morse-Smale assumption,
it follows that $\ker \Dd^\eps_u$
is generated by a nonzero vector
$Z^\eps$.
Because $\im {\Dd^\eps_u}^*$
can be identified with
the orthogonal complement
of $\ker \Dd^\eps_u$ with respect to
the $L^2$ inner product
$\langle\cdot,\cdot\rangle_\eps$,
it remains to prove that the function
$$
     \theta^\eps(\sigma)
    :=-\left\langle Z^\eps_\sigma,\zeta_\sigma
     \right\rangle_\eps,\qquad
     Z^\eps
    :=\begin{pmatrix}\p_su\\
     \Nabla{t}\p_su\end{pmatrix}
    -{\Dd^\eps_u}^*
     \left(\Dd^\eps_u{\Dd^\eps_u}^*\right)^{-1}
     \Dd^\eps_u
     \begin{pmatrix}\p_su\\\Nabla{t}\p_su
     \end{pmatrix}
$$
admits a zero.
Here the $s$-shift is defined by
$\zeta_\sigma(\cdot,\cdot)
:=\zeta(\cdot+\sigma,\cdot)$.

In Step four we prove surjectivity.
Assume by contradiction that there
is a sequence $(u_i,v_i)\in\Mm^{\eps_i}$
with $\eps_i$ converging to zero and such
that $(u_i,v_i)\notin\Tt^{\eps_i}(\Mm^0)$.
Viewing the $u_i$ as approximate zeroes
of the parabolic section $\Ff_0$, defined
by the left hand side of~(\ref{eq:heat}),
we construct a sequence of parabolic
solutions $u_i^0\in\Mm^0$ by
Newton-Picard iteration.
Here two conditions need to be satisfied
uniformly for large $i\in\N$.
Firstly, we need
a small initial value of $\Ff_0$.
This follows from the elliptic
equations~(\ref{eq:floer-eps})
and Step one:
$$
     \Norm{\Ff_0(u_i)}_p
     =\Norm{\p_su_i-\Nabla{t}\p_tu_i
     -\grad\Vv(u_i)}_p
     =\eps_i^2\Norm{\Nabla{t}\Nabla{s} v_i}_p
     \le c\eps_i^2.
$$
Secondly, one needs to prove
asymptotic decay of the form
$$
     \Abs{\p_su_i(s,t)}+\Abs{\Nabla{s}v_i(s,t)}
     \le\frac{c}{1+s^2}.
$$
This follows from uniform exponential decay
proved in Step two. The iteration shows
that $(u_i^0,\p_tu_i^0)$ and
the given elliptic solution
$(u_i,v_i)$ are sufficiently close
(for large $i$) such that the time shift
argument of Step three applies.
Hence $(u_i,v_i)=\Tt^{\eps_i}(u_i^0)$, but
this contradicts the assumption.

\section{Finite dimensional approximation}
\label{sec:V}

Reviewing Viterbo's
paper~\cite{V96}
is somewhat delicate
due to its state
of presentation.
As a way out we decided to
enlist the steps of proof
as provided by~\cite{V96}
using the original notation
and conventions.
In view of its independent
interest, we discuss
Step~1 including full details
(up to the hypothesis that
both gradings coincide)
in a separate section.

Throughout let $(M,g)$ be
a closed Riemannian manifold
of dimension $n$.
In~\cite{V96} cohomology
is considered and the main
result is stated in the form
\begin{equation}\label{eq:viterbo-main}
     \HF^*(DT^*M)
     \simeq\Ho^*(\Ll M)
\end{equation}
where $DT^*M$ denotes the
open unit disc bundle.
On both sides
contractible loops are
considered and the isomorphism
is claimed with rational
coefficients. A further claim is
a version of~(\ref{eq:viterbo-main})
for $S^1$-equivariant cohomologies
(see~\cite{V99} for applications).
Since orientation
of moduli space
is not discussed in~\cite{V96},
we assume
throughout Section~\ref{sec:V}
that all homologies take
coefficients in $\Z_2$.

Actually~\cite{V96} is part two
of a series of two papers.
The left hand side
of~(\ref{eq:viterbo-main})
is defined in part one~\cite{V99}
in the more general context
of symplectic manifolds with
contact type boundary (see the
excellent recent survey by
Oancea~\cite{O04}):
take the symplectic
completion of the closed
unit disc bundle
(which is symplectomorphic
to $(T^*M,\omega_0)$ itself),
fix $\delta<1$ close to $1$
and consider Hamiltonians of the form
$$
     H(t,q,p)=h_\lambda(\abs{p}).
$$
Here $h_\lambda$ is a smooth
convex real function
which vanishes for $\abs{p}<\delta$
and which is linear
of slope $\lambda>0$
for $\abs{p}\ge1$.
To define Floer cohomology
associated to $h_\lambda$
it is important to assume
that $\lambda$ is not
the length of a contractible
periodic geodesic in $M$
(see Remark~\ref{rmk:radial}).
The left hand
side of~(\ref{eq:viterbo-main})
is then defined
by the direct limit
$$
     \HF^*(DT^*M)
     :=\underset{\lambda\to\infty}
     {\underrightarrow{\lim}}
     \HF^*(T^*M,h_\lambda,J_g).
$$
To prove~(\ref{eq:viterbo-main})
it suffices to show
$$
     \HF^*(T^*M,h_\lambda,J_g)
     \simeq
     \Ho^*(\Ll^{\lambda^2/2} M).
$$
The proof rests on the idea
of Chaperon~\cite{Ch84}
to adapt the finite
dimensional approximation
of the loop space
via piecewise geodesics
to the case of Hamiltonian
flow lines. Givental~\cite{G89}
refined this idea for
Hamiltonian flows
which are periodic in time.

From now on fix
$h=h_\lambda$, denote
by $\varphi_t^h$ the time-t-map
generated by the Hamiltonian
vector field $X_h$, and
set $\varphi^h:=\varphi_1^h$.
The proof has seven steps:
\begin{equation*}
\begin{split}
     \HF^*(T^*M,h,J_g)
    &\stackrel{1}{\simeq}
     \HF^*(\Delta_r,\Gamma_r
     (\varphi^h),
     (T^*M\times
     \overline{T^*M})^{\times r},J_r
     ) \\
    &\stackrel{2}{\simeq}
     \HF^*(U_{r,\eps},
     \graph\: dS_\psi,
     T^*U_{r,\eps},\Psi_*J_r
     ) \\
\end{split}
\end{equation*}
\begin{equation*}
\begin{split}
    &\stackrel{3}{\simeq}
     \HM^{*+rn}(U_{r,\eps},
     S_\psi)\\
    &\stackrel{4}{\simeq}
     \HI^{*+rn}(U_{r,\eps},
     \nabla S_\psi)\\
    &\stackrel{5}{\simeq}
     \HI^{*+rn}(U_{r,\eps},
     \xi_\psi)\\
    &\stackrel{6}{\simeq}
     \Ho^*(\Lambda^a_{r,\eps}) \\
    &\stackrel{7}{\simeq}
     \Ho^*(\Ll^{\lambda^2/2} M).
\end{split}
\end{equation*}
Here $\HF^*(L_0,L_1;N)$ denotes
Floer cohomology associated to
Lagrangian submanifolds
$L_0$ and $L_1$
of a symplectic manifold $N$,
Morse cohomology
is denoted by $\HM^*$,
and the cohomological
Conley index by $\HI^*$.
The other symbols are
introduced as they appear
in our discussion of
the seven steps.
The whole proof relies
on writing $\varphi^h$
as an $r$-fold composition
$\psi^r:=\psi\circ\dots \circ\psi$
for sufficiently large $r\in\N$.
Since $h$ does
not depend on time,
we can indeed choose
$\psi:=\varphi_{1/r}^h$.

\smallbreak
\noindent
{\sc Step 1}
Let $\omega_0$ denote the
canonical symplectic structure
on $T^*M$, fix $r\in\N$, and consider
the symplectic manifold
$$
     (T^*M\times\overline{T^*M})^{\times r}
   :=
     \Bigl((T^*M\times T^*M)^{\times r},
     \omega_r:=\bigoplus_1^r
     \omega_0 \oplus -\omega_0\Bigr)
$$
and the Lagrangian submanifolds
$$
     \Gamma_r(\varphi^h)
   :=\{(z_1,\psi z_r;
        z_2,\psi z_1; \dots ;
        z_r,\psi z_{r-1})\mid
        z_1,\dots,z_r \in T^*M\},\quad
     \Delta_r
   :=\Delta^{\times r}.
$$
Here $\Delta$ denotes the diagonal
in $T^*M\times T^*M$.
Let $\Pp(h)$ denote
the set of contractible 1-periodic
orbits of $h$.
Since
$$
     \Gamma_r(\varphi^h)
     \cap\Delta_r
    =\{(z,z;\psi z,\psi z; \dots ;
     \psi^{r-1} z,\psi^{r-1} z) \mid
     z\in\Fix \:\varphi^h\}
     \simeq \Pp(h),
$$
it follows that the generators
of both chain groups
coincide up to natural identification.
That this identification
preserves the grading
given by the Conley-Zehnder index
and the Maslov index, respectively,
seems to be an open problem.\\
The boundary operator
of Lagrangian intersection
Floer homology counts $r$-tuples
of pairs of $J_r$-holomorphic strips
in $T^*M$, where $J_r$
is defined by~(\ref{eq:J_r}).
The boundary operator
of periodic
Floer homology counts
$J_g$-holomorphic cylinders in $T^*M$.
The idea to prove that both boundary
operators coincide is to establish
a one-to-one correspondence
by gluing together the strips
of an $r$-tuple to obtain a cylinder.
Details of this argument will be
discussed in a separate
section below.

\begin{remark}[Transversality]\rm
Working with time independent
Hamiltonians seems unrealistic
at a first glance, 
since Morse-Smale
transversality requires
time dependence.
We suggest perturbing $h$
in $\Cinf(\R/r^{-1}\Z\times DT^*M)$.
Due to the $1/r$-periodicity of
the perturbed Hamiltonian
$H$ the crucial decomposition
$\varphi^H=\psi^r$ is still available:
fix an initial time, say $t_0=0$,
and set $\psi:=\varphi^H_{0,1/r}$.
Since $H$ is still radial
outside a compact set,
Step~1 still goes through.
\end{remark}

\begin{figure}[ht]
  \centering
   \epsfig{figure=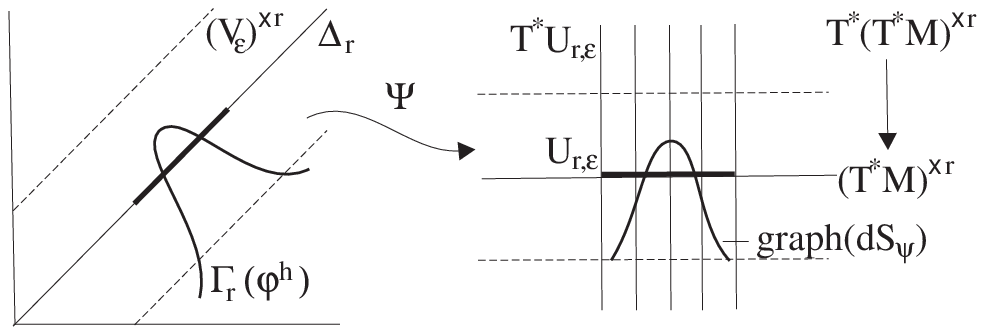}
   \caption{Change of coordinates $\Psi$
            identifying
            $\Gamma_r(\varphi^h)
            \cap(V_\eps)^{\times r}$
            with $\graph\:dS_\psi$.}
   \label{fig:fig-spsi}
\end{figure}

\smallbreak
\noindent
{\sc Step 2}
In this step coordinates
are changed via
a symplectomorphism which
identifies $\Delta_r$ with
$(T^*M)^{\times r}$ and
the essential part of
$\Gamma_r(\varphi^h)$ with the
graph of the differential of a function.
Let $d$ denote
the Riemannian distance on $M$ and
consider the neighborhood
$(V_\eps)^{\times r}$ of $\Delta_r$ where
$$
     V_\eps:=\{ (q,p;Q,P)\in
     T^*M\times T^*M \mid
     d(q,Q)\le \eps \}.
$$
Since $\Delta$ is diffeomorphic
to $T^*M$, clearly
$\Delta_r$ is diffeomorphic
to $(T^*M)^{\times r}$.
For sufficiently small $\eps>0$
a proper symplectic embedding
$$
     \Psi:(V_\eps)^{\times r}\to
     T^*(T^*M)^{\times r}
$$
is constructed
in~\cite[Lemma~1.1]{V97}
identifying $\Delta_r$ with
the zero section
(see Figure~\ref{fig:fig-spsi}).
According to~\cite{V97}
there exists $\eps>0$ sufficiently small
such that for every 
sufficiently large $r\in\N$
there is a function $S_\psi$
(denoted by $\Ss_\Phi$ in~\cite{V97};
note that $\psi$ depends on $r$)
which is defined on the set
$$
     U_{r,\eps}
     :=\{ (q_1,p_1;\dots;q_r,p_r)\in
     (T^*M)^{\times r} \mid
     d(q_j,q_{j+1})\le \frac{\eps}{2}\,\,
     \forall j\in\Z_r \}.
$$
It has the property that
the image under $\Psi$ of
the part of $\Gamma_r(\varphi^h)$
contained in $(V_\eps)^{\times r}$
equals the graph of $dS_\psi$.

\noindent
To see existence of $S_\psi$
note that the Lagrangian
$\Psi\left(\Gamma_r(\id)
\cap(V_\eps)^{\times r}\right)$
is a graph over $U_{r,\eps}$ and
$\Psi\left(\Gamma_r(\varphi^h)
\cap(V_\eps)^{\times r}\right)$
is close to this graph
whenever $r$ is sufficiently large.
Hence it is a Lagrangian
graph over $U_{r,\eps}$ itself.
It is exact,
because the symplectomorphism
$\varphi^h$ is Hamiltonian.
To arrive at the claim of Step~2
consider the two
isomorphisms given by
\begin{equation*}
\begin{split}
    &\HF^*(\Delta_r,\Gamma_r
     (\varphi^h),
     (T^*M\times
     \overline{T^*M})^{\times r},J_r
     ) \\
    &\simeq
     \HF^*((T^*M)^{\times r},
     \Psi(\Gamma_r(\varphi^h)
     \cap(V_\eps)^{\times r}),
     T^*(T^*M)^{\times r},\Psi_*J_r
     ) \\
    &\simeq
     \HF^*(U_{r,\eps},
     \graph\: dS_\psi,
     T^*U_{r,\eps},\Psi_*J_r
     ).
\end{split}
\end{equation*}
The first isomorphism
is induced by the symplectic
embedding $\Psi$.
One needs to prove
that the Floer complex
associated to
$(\Delta_r,\Gamma_r(\varphi^h);
(T^*M\times\overline{T^*M})^{\times r};
J_r)$ lives entirely inside
the neighborhood
$(V_\eps)^{\times r}$ of $\Delta_r$.
This is clear for the
generators, i.e. the
points of intersection,
but not so much for
the connecting trajectories,
i.e. the $J_r$-holomorphic strips.
The second isomorphism
follows, if again the whole Floer complex
associated to the large
space $T^*(T^*M)^{\times r}$
lives in the smaller
space $T^*U_{r,\eps}$.
This is clear for the
generators and for
the connecting trajectories
a proof is given in~\cite{V96}.

\smallbreak
\noindent
{\sc Step 3}
The Hamiltonian
flow on $T^*U_{r,\eps}$
generated by the Hamiltonian
$(X,Y)\mapsto -S_\psi(X)$
is given by
$$
     \phi_t(X,Y)=(X,Y+tdS_\psi(X)).
$$
According to~\cite{V96}
the Floer cohomologies
associated to $\Psi_*J_r$
and $(\phi_t)_*\Psi_*J_r$,
respectively, are isomorphic
by continuation.
Furthermore, it is
shown in~\cite{V96}
that the generators and
connecting trajectories
of the Floer complex
associated to $(\phi_t)_*\Psi_*J_r$
are in one-to-one correspondence
with those of the Morse complex
associated to $(U_{r,\eps},\nabla S_\psi)$.
Therefore the corresponding
cohomologies coincide up to
the shift in the grading.

\smallbreak
\noindent
{\sc Step 4}
Here $\HI^*$ denotes the
cohomological Conley index 
and the isomorphism is refered
to~\cite{F89b}.

\smallbreak
\noindent
{\sc Step 5}
Existence of the pseudogradient vector
field $\xi_\psi$ for $S_\psi$
is established in~\cite{V97}.
Let ${\rm I}_*(U_{r,\eps},\xi_\psi)$ be
the Conley index, i.e.
the homotopy type
of the quotient of $U_{r,\eps}$
by the exit set with respect
to the flow generated by $\xi_\psi$.
According to~\cite{V96}
it holds
$$
     {\rm I}_*(U_{r,\eps},\nabla S_\psi)
     \simeq
     {\rm I}_*(U_{r,\eps/2},\xi_\psi)
$$
and the spaces are independent
of $\eps$ whenever $r$
is sufficiently large.

\smallbreak
\noindent
{\sc Step 6}
Let $\eps>0$ be smaller than the injectivity
radius of $M$ and consider the finite
dimensional approximation
of the free loop space given by
$$
     \Lambda_{r,\eps}
     :=\{ (q_1,\dots,q_r)\in
     M^{\times r} \mid
     d(q_j,q_{j+1})\le \frac{\eps}{2}\,\,
     \forall j\in\Z_r \}.
$$
Define a function on $\Lambda_{r,\eps}$ by
$$
     E_\psi(q_1,\dots,q_r)
     :=
     \sup_{(p_1,\dots,p_r)}
     S_\psi(q_1,p_1;\dots;q_r,p_r)
$$
and set
$\Lambda_{r,\eps}^a
:=\{E_\psi\le a\}$.
The Conley index of
$(U_{r,\eps},\xi_\psi)$
is calculated 
in~\cite[Prop.~1.7]{V97}
and the result is
the Thom space of some
vector bundle of rank $rn$
over $\Lambda_{r,\eps}^a$,
for sufficiently large $a$.
Taking cohomology
the Thom isomorphism
leads to the shift 
in the grading and proves
Step~6.

\smallbreak
\noindent
{\sc Step 7}
According to~\cite[p.~438]{V97}
the space $\Lambda_{r,\eps}^a$ approximates
the free loop space
$\Ll^{\lambda^2/2} M$ for $a\to \infty$,
$r\to \infty$ and $\eps\to0$.

\smallbreak
Note that
generating function homology
has not been used throughout
the seven steps -- in contrast
to what one expects after
a glimpse into~\cite{V96}.

\subsection*{Step~1 revisited}

We shall consider the
case of Floer homology
instead of Floer cohomology
and continue our discussion of
Step~1 above.
First we review
the original
approach in~\cite{V96}.
Then we propose a
proof of the final argument
along different lines.

\begin{figure}[ht]
  \centering
   \epsfig{figure=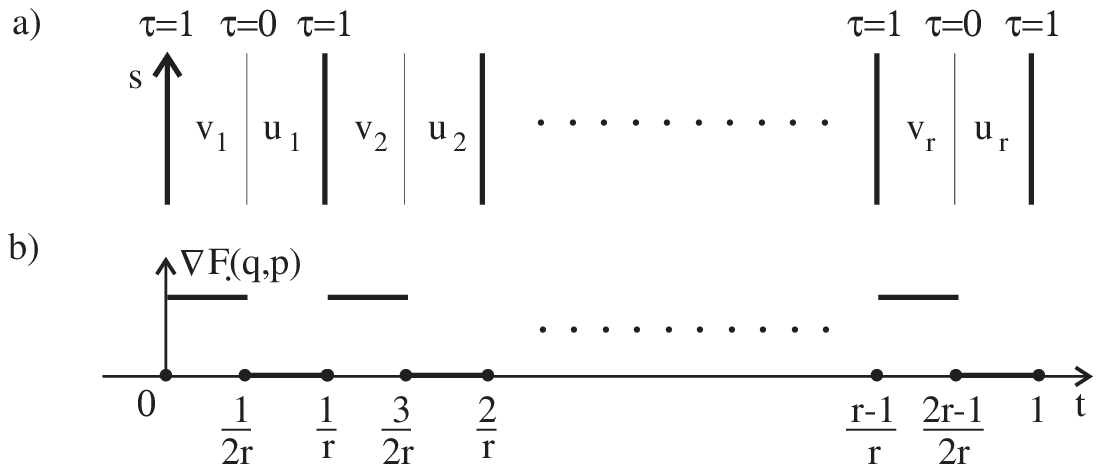}
   \caption{a) Cylinder of
            $J$-holomorphic strips.
            b) Time dependence of $\nabla F$.}
   \label{fig:fig-wjuj}
\end{figure}

To establish the isomorphism
between periodic and Lagrangian
Floer homology
it remains to check
that the connecting trajectories
arising in both situations
are in one-to-one correspondence.
In the case of periodic Floer homology
these are solutions
$w:\R\times S^1\to T^*M$
satisfying
\begin{equation}\label{eq:floer-H}
     \p_s w+J_g(w)\p_t w
     -\nabla h(w)
     =0
\end{equation}
and appropriate
boundary conditions.
In the case of Lagrangian intersection
Floer homology these
are $r$ pairs of strips
$
     W:=(u_1,\hat{v}_1;\dots;
     u_r,\hat{v}_r) : 
     \R\times [0,1]\to 
     (T^*M\times\overline{T^*M})^{\times r}
$
satisfying
\begin{equation}\label{eq:U-V}
     \p_s W+J_r(W)
     \p_\tau W
     =0,\qquad
     W(s,0)\in\Delta_r,\qquad
     W(s,1)\in
     \Gamma_r(\varphi^h).
\end{equation}
Here we suggest to use the
almost complex structure
\begin{equation}\label{eq:J_r}
     J_r=J_{r,\tau}:=\bigoplus_1^r J_g\oplus 
     -(\varphi_{-\frac{\tau}{r}}
     ^h)^* J_g
\end{equation}
in order to make the key
idea in~\cite{V96} work:
achieve matching boundary
conditions by redefining the
$\hat{v}_j$'s
in the form
$v_j(s,\tau):=
\varphi_{-\frac{\tau}{r}}^h\circ
\hat{v}_j(s,\tau)$.
With this definition~(\ref{eq:U-V})
is equivalent to
\begin{equation}\label{eq:u-v}
\begin{aligned}
     \p_s u_j+J_g(u_j)\p_\tau u_j
    &=0,&
     v_j(s,0)
    &=u_j(s,0),
    \\
     \p_s v_j-J_g(v_j)\p_\tau v_j
    &=\frac{1}{r}\nabla 
     h(v_j),\qquad&
     u_j(s,1)
    &=v_{j+1}(s,1),
\end{aligned}
\end{equation}
for $j=1,\dots,r$.
Here and throughout we
identify $r+1$ and $1$.
Note that the minus sign in the second
PDE is fine, since we need to
reverse time in the $v_j$'s
when fitting them together
with the $u_j$'s to obtain
a perturbed $J$-holomorphic map
$\tilde{w}:\R\times S^1\to T^*M$
(see Figure~\ref{fig:fig-wjuj}a).
More precisely, define
\begin{equation}\label{eq:w-tilde}
      \tilde{w}(s,t)
      :=
     \begin{cases}
      v_j(2rs,1-2rt+2j-2)
     &,\, j\in\{1,\dots,r\},\,
      t\in [\frac{2j-2}{2r},
      \frac{2j-1}{2r}],
     \\
      u_j(2rs,2rt-2j+1)
     &,\, j\in\{1,\dots,r\},\,
      t\in [\frac{2j-1}{2r},
      \frac{2j}{2r}],
     \end{cases}
\end{equation}
and consider the perturbation
associated to the Hamiltonian
(see Figure~\ref{fig:fig-wjuj}b)
\begin{equation*}
      F_t(q,p)
      :=
     \begin{cases}
      2h(\abs{p})
     &,\, j\in\{1,\dots,r\},\,
      t\in [\frac{2j-2}{2r},
      \frac{2j-1}{2r}],
     \\
      0
     &,\,\text{else}.
     \end{cases}
\end{equation*}
The argument
in~\cite{V96} concludes as follows:
the time-1-maps associated to
$F$ and $h$ coincide, hence
$\Pp(F)\simeq \Pp(h)$,
and $\tilde{w}$ solves
\begin{equation*}
     \p_s \tilde{w}
     +J_g(\tilde{w})\p_t \tilde{w}
     -\nabla F_t(\tilde{w})
     =0
\end{equation*}
iff the pairs
$(u_j,v_j)$ solve~(\ref{eq:u-v}).
Hence both
chain complexes are equal and
$$
     \HF_*(\Delta_r,
     \Gamma_r(\varphi^h),
     (T^*M\times
     \overline{T^*M})^{\times r},J_r)
     \simeq
     \HF_*(T^*M,F,J_g).
$$
According to~\cite{V96}
continuation\footnote{After writing
this paper Viterbo informed us 
that here continuation
does not refer to Floer continuation
but to the following:
assume the 1-periodic orbits do not depend 
on the parameter $\lambda\in[0,1]$
of a homotopy $f_\lambda$
between $F$ and $h$ and
consider two 1-periodic orbits $z^\pm$
of index difference one. Set
$X=\{(\lambda,w)\mid\p_s w
+J_g(w)\p_t w
=\nabla f_\lambda(w),\;
\lim_{s\to \pm \infty} w(s,t)
=z^\pm(t)\}$.
Let $\R$ act on $X$ by $s$-shift of $w$.
Then the projection
$\pi:X/\R\to[0,1]$,
$(\lambda,[w])\mapsto\lambda$,
is Fredholm of index $0$.
Hence the algebraic numbers of
$\pi^{-1}(0)$ and $\pi^{-1}(1)$ are equal,
provided we have 
compactness control.
But these numbers represent
the boundary matrix elements between
$z^-$ and $z^+$ in
the Floer complexes associated
to $F$ and $h$, respectively.
In the general case, as long as
the 1-periodic orbits do not bifurcate, 
the same argument applies: 
the boundary conditions now 
depend on $\lambda$, but 
remain nondegenerate of constant 
index for all $\lambda$.
} shows
$
     \HF_*(T^*M,F,J_g)
     \simeq
     \HF_*(T^*M,h,J_g)
$.

\medskip
The crucial point
in defining the Floer homology
associated to $F$
and constructing a Floer
continuation homomorphism is
to prove apriori $C^0$ estimates
for the corresponding solutions.
Both problems are
nonstandard in the sense that
$F$ as well as
the homotopy of Hamiltonians
are \emph{asymptotically
nonconstant}.
Whereas a convexity
argument can probably
be adapted to solve
the first problem, we don't
see how to get the
$C^0$ estimates needed
for Floer continuation.
Hence we propose a proof along
different lines:
in Remark~\ref{rmk:nonsmooth}
we switch on the Hamiltonian
perturbation smoothly --
simply to remain
in the familiar
setting of \emph{smooth}
Hamiltonians.
Then, in Remark~\ref{re:continuation},
we avoid asymptotically
nonconstant Floer continuation
altogether by introducing
an intermediate step. All
$C^0$ estimates we need will follow
from~\cite[Proposition~2.3]{JOA4}.
For convenience we
recall this result below.
Roughly speaking, it asserts
that a Floer cylinder
whose ends are located
\emph{inside} $D T^*M$
cannot leave this set,
whenever outside the Hamiltonian is
radial with second
derivative bounded below or above.

\begin{proposition}[\cite{JOA4}]
\label{pr:subsolution}
Given $R,c\ge0$,
let $f\in\Cinf(\R\times S^1\times [R,\infty))$
satisfy $\p_s f_{s,t}'\ge0$ and
$f_{s,t}''\ge -c$ (or $f_{s,t}''\le c$)
for all $s$ and $t$.
Here $f_{s,t}(r):=f(s,t,r)$
and $f_{s,t}':=\frac{d}{dr} f_{s,t}$.
Assume further that
$H\in C^\infty(\R\times S^1\times TM)$
satisfies
$H(s,t,x,y)=f_{s,t}(\Abs{y})$
whenever $\Abs{y}\ge R$,
that the pair
$(u,v)\in C^\infty(\R\times S^1,TM)$
satisfies
\begin{equation*}\label{eq:l2}
     \begin{pmatrix}
       \p_su-\Nabla{t}v \\
       \Nabla{s}v+\p_tu
     \end{pmatrix}
     -\nabla H(s,t,u,v) =0,
\end{equation*}
and that there exists $T>0$ such that
$\Abs{v(s,\cdot)}\le R$
whenever $\Abs{s}\ge T$.
Then $\Abs{v}\le R$
on $\R\times S^1$.
\end{proposition}

\begin{remark}[Smooth Hamiltonian]
\label{rmk:nonsmooth}\rm
We smoothly switch off the Hamiltonian
perturbation near the boundary
of the $\hat{v}_j$-strips:
fix a nondecreasing smooth map
$\alpha:[0,1]\to\R$ with
$\alpha\equiv0$ near 0 and
$\alpha\equiv\frac{1}{r}$ near 1
(see Figure~\ref{fig:fig-alph}).
\begin{figure}[ht]
  \centering
   \epsfig{figure=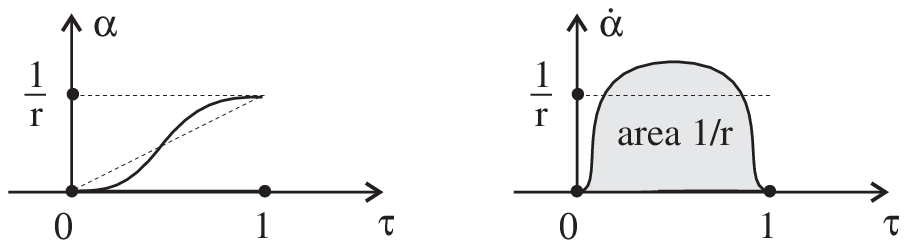}
   \caption{The function $\alpha$
            and its derivative.}
   \label{fig:fig-alph}
\end{figure}

\noindent
Define
$$
     v_j(s,\tau)
     :=\varphi^h
     _{-\alpha(\tau)}\circ
     \hat{v}_j(s,\tau),\qquad
     J_r=J_{r,\tau}:=\bigoplus_1^r J_g\oplus 
     -(\varphi^h
     _{-\alpha(\tau)})^* J_g.
$$
Then the Lagrangian boundary value
problem~(\ref{eq:U-V})
is equivalent to
\begin{equation}\label{eq:u-v-2}
\begin{aligned}
     \p_s u_j+J_g(u_j)\p_\tau u_j
    &=0,&
     v_j(s,0)
    &=u_j(s,0),
    \\
     \p_s v_j-J_g(v_j)\p_\tau v_j
    &=\dot\alpha(\tau) \nabla 
     h(v_j),\qquad&
     u_j(s,1)
    &=v_{j+1}(s,1),
\end{aligned}
\end{equation}
for $j=1,\dots,r$.
Define a $1/r$-periodic
function
(see Figure~\ref{fig:fig-fnew})
by
\begin{equation*}
      \beta(t)
      :=
     \begin{cases}
      2r\dot\alpha(1-2rt+2j-2)
     &,\, j\in\{1,\dots,r\},\,
      t\in [\frac{2j-2}{2r},
      \frac{2j-1}{2r}],
     \\
      0
     &,\,\text{else}.
     \end{cases}
\end{equation*}
\begin{figure}[ht]
  \centering
   \epsfig{figure=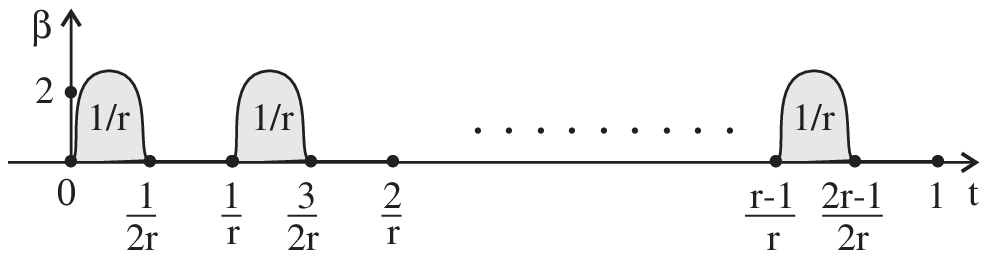}
   \caption{Smooth time dependence
            of Hamiltonian 
            perturbation.}
   \label{fig:fig-fnew}
\end{figure}

\noindent
Consider
the Hamiltonian $(t,q,p)\mapsto
\beta(t) h(\abs{p})$, and
let $\tilde{w}:\R\times S^1\to T^*M$
be given by~(\ref{eq:w-tilde}).
Then~(\ref{eq:u-v-2}) is equivalent to
$$
     \p_s \tilde{w}+J_g(\tilde{w})\p_t \tilde{w}
     -\beta(t) \nabla h(\tilde{w}) 
     =0.
$$
With these definitions
the Hamiltonian perturbations
on the strips $(u_j,v_j)$
fit together smoothly. Indeed
$\nabla (\beta h)=\beta\nabla h$
depends smoothly on $t$
(see Figure~\ref{fig:fig-fnew}).
\\
The fact that
$\int_0^{1/r}\beta(t)\: dt=1/r$ and
the identity $X_{\beta h}=\beta X_h$
together imply that the time-$1/r$-maps
associated to $\beta h$
and $h$, respectively, are equal.
This proves
$\Fix\:\varphi^h=
\Fix\:\varphi^{\beta h}$ and
$\Pp(h)\simeq\Pp(\beta h)$.
The latter correspondence is given by
mapping $z\in\Pp(h)$ to
$\tilde{z}(\cdot):=z(\int_0^\cdot \beta(t)\:dt)$.
Hence both chain complexes are equal
and therefore
$$
     \HF_*(\Delta_r,
     \Gamma_r(\varphi^h),
     (T^*M\times
     \overline{T^*M})^{\times r},
     J_r;\Z_2)
     \simeq
     \HF_*(T^*M,\beta h,J_g;\Z_2).
$$
Again we remark that
it is an open problem to
prove equality of the gradings.
Since the Hamiltonian $\beta h$
is not well behaved at
infinity, definition of the
right hand side
requires an additional
argument to obtain a
uniform $C^0$-bound for
Floer trajectories:
given our knowledge that all elements
of $\Pp(\beta h)$ take values in
$DT^*M$, the $C^0$ estimate
for radial time-dependent Hamiltonians
Proposition~\ref{pr:subsolution}
shows that \emph{all} Floer cylinders
connecting elements of $\Pp(\beta h)$
also take values in 
the bounded set $DT^*M$.
\end{remark}

\begin{remark}[Floer continuation]
\label{re:continuation}\rm
To conclude the proof of Step~1
we need to show that the Floer
homologies of $\beta h$ and $h$,
respectively, are isomorphic.

Choosing a different function
$\beta$, if necessary,
we may assume without loss of
generality that 
$\beta_m:=\norm{\beta}_\infty=3$.
(Any real number strictly larger
than $2$ can be realized
as such a maximum.)
Pick a sufficiently large real
$R_0>1$ and a smooth nondecreasing
cutoff function $\rho$
which equals 0 on $(-\infty,1]$
and 1 on $[R_0,\infty)$
such that the radial Hamiltonian
$$
     h^{\rho,\beta(t)}(r)
     :=\Bigl( \rho(r)+\beta(t)
     \bigl(1-\rho(r)\bigr)\Bigr)
     h(r),\qquad r=\abs{p},
$$
(see Figure~\ref{fig:fig-hrho})
is nondecreasing for every $t\in[0,1]$.
To check that such $R_0$ and $\rho$
exist is left as an exercise.
(Hint: start with
$$
     R_0:=-11\frac{c_\lambda}{\lambda}+12,
     \qquad
     \tilde{\rho}(r)
     :=\begin{cases}
     0
    &, r\le 1 \\
     -\frac{(r-R_0)^2}{(1-R_0)^2}+1
    &, r\in [1,R_0] \\
     1
    &, r\ge R_0
     \end{cases}
$$
and check that indeed
$(h^{\tilde{\rho},\beta(t)})^\prime\ge 0$
for every $t\in[0,1]$.
Then smooth out
$\tilde{\rho}$ near $r=1$ and $r=R_0$
using cutoff functions
whose derivatives are
supported in $[1,1+\mu]$
and $[R_0-\mu,R_0]$, respectively,
for appropriate $\mu>0$.)
In the $t$-variable
$h^{\rho,\beta(t)}$ oscillates
between $h^{\rho,\beta_m}$
and $\rho h$
(see Figure~\ref{fig:fig-hrho}).
For $r\le 1$ the Hamiltonian
$h^{\rho,\beta}$
coincides with $\beta h$
and for $r\ge R_0$
with $h(r)=-c_\lambda+\lambda r$.
\begin{figure}[ht]
  \centering
   \epsfig{figure=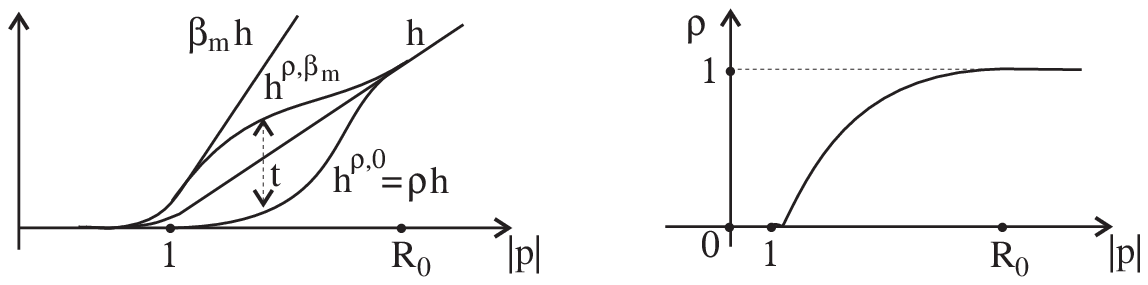}
   \caption{The family of Hamiltonians
            $h^{\rho,\beta(t)}$
            and a cutoff function $\rho$.}
   \label{fig:fig-hrho}
\end{figure}

Since $h$ and $h^{\rho,\beta}$ are both
linear of slope $\lambda$
outside $D_{R_0}T^*M$,
any homotopy $f_s$
given by a convex combination
of the two Hamiltonians satisfies
$f_s^\prime\equiv \lambda$
outside $D_{R_0}T^*M$.
Hence the two conditions
$f_s^{\prime\prime}\equiv 0$
and $\p_s f_s^\prime\equiv 0$
in Proposition~\ref{pr:subsolution}
are satisfied for $R=R_0$
and the homotopy
gives rise to a Floer continuation
homomorphism. It is an isomorphism
by the reverse homotopy argument.
Hence
$$
     \HF_*(T^*M,h,J_g;\Z_2)
     \simeq
     \HF_*(T^*M,h^{\rho,\beta},J_g;\Z_2).
$$

By the end of Remark~\ref{rmk:nonsmooth}
all elements of $\Pp(\beta h)$
and all connecting Floer
trajectories associated to $\beta h$
take values in $DT^*M$.
We claim that all elements
of $\Pp(h^{\rho,\beta})$
take values
in $DT^*M$, too.
Then by Proposition~\ref{pr:subsolution}
the same will be true
for all connecting
Floer cylinders associated to
$h^{\rho,\beta}$.
But on $S^1\times DT^*M$
both Hamiltonians coincide,
hence both chain complexes are equal
and
$$
     \HF_*(T^*M,h^{\rho,\beta},J_g;\Z_2)
     \simeq
     \HF_*(T^*M,\beta h,J_g;\Z_2).
$$

To prove the claim
consider the Hamiltonian vector field
$X_{\abs{p}^2/2}$ which
generates the geodesic flow on $T^*M$.
Since the Hamiltonian $h^{\rho,\beta}$
is radial, it holds
$$
     X_{h^{\rho,\beta}}(q,p)
     =r^{-1} (h^{\rho,\beta})^\prime (r)
     X_{r^2/2}(q,p),\qquad r=\abs{p}.
$$
This shows that the projection
to the zero section $M$
of a $X_{h^{\rho,\beta}}$-trajectory
is a reparametrized geodesic.
In particular, if
$(x,y)\in\Pp(h^{\rho,\beta})$,
then $x:S^1\to M$ is a closed geodesic
(not necessarily parametrized
with respect to arc length)
and $\abs{y(t)}$ is
independent of $t$. Assume
by contradiction
$r:=\abs{y(t)}\ge 1$, then the
length of the closed geodesic $x$
is given by
\begin{equation*}
\begin{split}
     \ell(x)
    &=\int_0^1\Abs{\dot x(t)} dt
     =\int_0^1\Abs{
     \frac{(h^{\rho,\beta})^\prime(r)}{r} 
     g(x(t))^{-1} y(t)} dt \\
    &=\int_0^1\Abs{
     (h^{\rho,\beta})^\prime(r)} dt
     =\int_0^1
     \rho^\prime(1-\beta)h
     +\bigl(\rho+\beta(1-\rho)\bigr)
     \lambda \,dt \\
    &=\rho^\prime h\int_0^1(1-\beta)\, dt
     +\rho\lambda\int_0^1(1-\beta)\, dt
     +\lambda \int_0^1\beta\, dt
     =\lambda.
\end{split}
\end{equation*}
Here we used $(h^{\rho,\beta})^\prime\ge 0$,
the fact that $h^\prime(r)=\lambda$
whenever $r\ge 1$, and $\int_0^1 \beta=1$.
But $\lambda$ was chosen in the
complement of the length spectrum of $M$.
\end{remark}

\begin{remark}[Closed aspherical case]\rm
Replacing the cotangent bundle
by a closed \emph{symplectically
aspherical} manifold $(N,\omega)$,
i.e. $\omega$ and the first
Chern class vanish over $\pi_2(N)$,
the argument of Step~1 in the
case $r=1$ provides a new
proof of the known fact
$$
     HF_*(N,H;\Z)
     \simeq
     HF_*(\Delta,\graph\:\varphi^H;
     N\times \overline{N};\Z).
$$
Note that
Remark~\ref{re:continuation}
can be replaced by
taking right away a homotopy
between $\beta H$ and $H$
and applying the standard
reverse homotopy continuation argument.
It is interesting to compare with
the proof in~\cite{BPS03}.
\end{remark}

\section{Example: The euclidean torus}
\label{sec:torus}

We compute $L^2$ Morse homology
of the loop space
of the euclidean torus $\T^n=\R^n/2\pi\Z^n$
and Floer homology of its cotangent bundle.
Since the differential equations decouple,
the calculation essentially reduces to
the case $n=1$.

\subsubsection*{\boldmath$L^2$ Morse homology 
                for $\Ll_\alpha\T^1$}
\label{subsubsec:HM(T^1)}

Given $\alpha\in\Z$,
we think of $x\in\Ll_\alpha\T^1$ as a smooth
map $x:\R\to\R$ satisfying
$x(t+1)=x(t)+2\pi\alpha$.
Hence $\alpha$
is the winding number of $x$.
We fix $\alpha$ and
analyze the problem
for each component of the
free loop space separately.

The equation for the free pendulum
is $\p_t\p_t x=0$
and its solutions are given by
$\{ x(t)=2\pi\alpha t+x_0\mid
x_0\in\T^1\}$.
Upon introducing a
time-dependent potential
on $S^1\times\T^1$ of the form
$
     V_\alpha(t,q)
     :=-\cos \left(q-2\pi\alpha t\right)
$,
the circle of solutions splits into two
nondegenerate critical points.
Equation~(\ref{eq:crit})
for the critical points of
the functional
$\Ss_{V_\alpha}$ defined on
$\Ll_\alpha \T^1$ becomes
\begin{equation}\label{eq:crit-1}
     \p_t\p_t x
     =-\sin\left(x-2\pi\alpha t\right).
\end{equation}
Observe that the case $\alpha=0$ describes
the mathematical pendulum with gravity.
Obvious 1-periodic solutions
of~(\ref{eq:crit-1}) are
the two equilibrium states
\emph{pendulum up} and
\emph{pendulum down}.
They are described by
the constant functions $x^{(1)}=\pi$
and $x^{(0)}=0$, respectively.
For general $\alpha$
equation~(\ref{eq:crit-1})
still describes
the same mathematical pendulum,
but now the observer rotates
with constant angular
speed $2\pi\alpha$.
The two equilibrium states 
of the pendulum are then
seen by the observer
as rotations and described by
\begin{equation}\label{eq:equilibrium}
     x^{(1)}=2\pi\alpha t+\pi,\qquad
     x^{(0)}=2\pi\alpha t.
\end{equation}
Of course, the pendulum
admits plenty of nonstationary
periodic solutions, but their
periods are strictly bigger
than one. Hence $\Ss_{V_\alpha}$
has only the two critical points
$x^{(1)}$ and $x^{(0)}$.
Their actions are
$2\pi^2\alpha^2+1$ and
$2\pi^2\alpha^2-1$, respectively.

To compute the Morse index of $x^{(1)}$
we linearize equation~(\ref{eq:crit-1})
at $x^{(1)}$ and solve the
resulting eigenvalue problem
$-\p_t\p_t \xi-\xi=\lambda\xi$
for $\xi:\R\to\R$ with
$\xi(t+1)=\xi(t)$.
The solution to this ODE is
given by $\xi(t)=\sin \sqrt{\lambda+1}t$
and the periodicity condition implies
$\lambda\in\{4\pi^2k^2-1\mid
k=0,1,2,\dots\}$.
This shows $\IND_{V_\alpha}(x^{(1)})=1$.
Linearizing at $x^{(0)}$ we arrive at
$-\p_t\p_t\xi+\xi=\lambda\xi$
with solution
$\xi(t)=\sin \sqrt{\lambda-1}t$
and $\lambda\in\{4\pi^2k^2+1\mid
k=0,1,2,\dots\}$. Hence
$\IND_{V_\alpha}(x^{(0)})=0$.
Therefore the only
nontrivial chain groups
are given by
\begin{equation}\label{eq:chain-1}
     C_0({V_\alpha})
    =\Z x^{(0)},\qquad
     C_1({V_\alpha})
    =\Z x^{(1)}.
\end{equation}
There must be precisely two connecting
trajectories up to $s$-shift,
because the unstable manifold
of $x^{(1)}$ is 1-dimensional.
For $\sigma\in C^\infty(\R,\R)$
consider the \emph{Ansatz}
$u(s,t):=2\pi\alpha t+\sigma(s)$.
Then the heat equation~(\ref{eq:heat})
is equivalent to
\begin{equation}\label{eq:sigma}
     \sigma^\prime(s)=-\sin \sigma(s).
\end{equation}
The initial conditions
$\sigma(0)=0$ and $\sigma(0)=\pi$
produce the constant trajectories
$u(s,t)=x^{(0)}(t)$
and $u(s,t)=x^{(1)}(t)$, respectively.
The two 1-parameter families
of connecting trajectories
are obtained by choosing
$\sigma(0)$ in $(0,\pi)$ or in
$(-\pi,0)$.
Choosing an orientation of the unstable
manifold of $x^{(1)}$ induces an
orientation of each of the
two 1-parameter families.
For one family this orientation
coincides with the flow direction
and for the other one it does not,
so one characteristic sign
is positive and one is negative.
Hence
$\p^M x^{(1)}=x^{(0)}-x^{(0)}=0$.
This proves that
$L^2$ Morse homology of
$(\Ll_\alpha \T^1,\Ss_{V_\alpha})$
is in fact given by its chain
groups~(\ref{eq:chain-1}).

\subsubsection*{Floer homology 
                of \boldmath$T^*\T^1$}
\label{subsubsec:HF(T^*T^1)}

Let $V_\alpha$ be as above
and define the Hamiltonian
$H_{V_\alpha}(t,q,p):=\frac12 p^2
+V_\alpha(t,q)$ for
$(t,q,p)\in S^1\times\T^1\times\R$.
Hamilton's equations and their
solutions are given by
$$
     \begin{pmatrix}\dot x(t)\\
     \dot y(t)\end{pmatrix}
    =\begin{pmatrix} y(t)\\
     -\sin\left( x(t)-2\pi\alpha t\right)
     \end{pmatrix},\;\;
    z^{(1)}(t)
    =\begin{pmatrix} 2\pi\alpha t+\pi\\
     2\pi\alpha\end{pmatrix},\;\;
    z^{(0)}(t)
    =\begin{pmatrix} 2\pi\alpha t\\
     2\pi\alpha\end{pmatrix}.
$$
To calculate the two connecting
trajectories let $u$ be as above
and set $v(s,t):=2\pi\alpha$.
Floer's equations~(\ref{eq:heat})
are then equivalent to~(\ref{eq:sigma})
and the argument continues as above
showing that
$\HF_*(T^*\T^1,H_{V_\alpha},J_1;\alpha)$
equals the chain
groups~(\ref{eq:chain-1}).

One might wonder what happens in the case of
the $\eps$-Floer equations~(\ref{eq:floer-eps}).
With our Ansatz for $u$ and $v$ these are
again equivalent to~(\ref{eq:sigma}),
in particular $\eps$ disappears.
This means that the map
$\Tt^\eps:\Mm^0\to\Mm^\eps$
is essentially the identity.
It maps $u$ to $(u,\p_tu)$
for any $\eps>0$.
The case $\eps=1$ and our
knowledge in the Morse case
show that up to $s$-shift
there are precisely
two Floer trajectories

\subsubsection*{Higher dimensional case}
\label{subsubsec:general}

Fix $n\ge1$ and
$\alpha=(\alpha_1,\dots,\alpha_n)\in\Z^n$.
On the euclidean torus
$\T^n=\R^n/\Z^n$ consider
the potential
$V_\alpha(t,q)$
given by the sum of $n$
one-dimensional potentials.
The critical points of
$\Ss_{V_\alpha}$ on
$\Ll_\alpha\T^n$
are smooth functions
$x:\R\to\R^n$ with
$x(t+1)=x(t)+2\pi\alpha$
satisfying the ODE
$$
     \p_t\p_t x(t)
     =-\nabla V_\alpha(t,x(t))
     =\left(
     -\sin\left(x_1(t)-2\pi\alpha_1 t\right)
     ,\dots,
     -\sin\left(x_n(t)-2\pi\alpha_n t\right)
     \right).
$$
Its solutions
are of the form
$(x^{(\nu_1)},\dots,
x^{(\nu_n)})$ where the components
$x^{(\nu_k)}$ are given by one of the
two equilibrium states
in~(\ref{eq:equilibrium}).
Hence the assignment
$(x^{(\nu_1)},\dots,
x^{(\nu_n)})\mapsto
(\nu_1,\dots,\nu_n)=:\nu$
identifies the set of critical points with
the direct sum $(\Z_2)^n$
of $n$ copies of $\Z_2$,
and the number of
critical points is $2^n$.
Moreover, the Morse index of
$x^{(\nu)}:=(x^{(\nu_1)},\dots,x^{(\nu_n)})$
equals $\Abs{\nu}$ and
\begin{equation}\label{eq:chain-n}
     C_k(V_\alpha)
    =\bigoplus_{\nu\in(\Z_2)^n,\,
     \Abs{\nu}=k}\;
     x^{(\nu)}\Z
    \simeq \Z^{\binom{n}{k}}.
\end{equation}
Arguing as above one observes that
trajectories  connecting critical
points of index difference one
come in pairs of opposite
characteristic signs
showing that all
boundary operators are zero.
(For example, in the case $n=3$
the critical point $(1,0,1)$
admits two connecting trajectories
to $(0,0,1)$ and two to $(1,0,0)$).
It follows that
$L^2$ Morse homology and
Floer homology both
coincide with the chain groups~(\ref{eq:chain-n})
and, by Theorem~\ref{thm:main},
so does $\Ho_*(\Ll_\alpha \T^n)$.

\section{Application:
         Noncontractible periodic orbits}
\label{sec:noncontractible}

The search for noncontractible
1-periodic orbits is a fairly recent
branch of symplectic geometry
pioneered by Gatien-Lalonde~\cite{GL00}
and Biran-Polterovich-Salamon~\cite{BPS03}.
(See also~\cite{L03} for some
generalizations of~\cite{GL00}).
These authors obtain existence results
under certain geometric assumptions,
for instance flatness of the metric.
More precisely, Theorem~\ref{thm:exist-orbit}
below is proved in~\cite{BPS03}
for the euclidean $n$-torus
and in the case of negative curvature
using the known simple structure
of the set of periodic geodesics.
In contrast our extension
to the case of \emph{arbitrary}
Riemannian metric
is based on Theorem~\ref{thm:main}.

Let $M$ be a closed connected
Riemannian manifold. Given a
homotopy class $\alpha$ of free
loops in $M$, define
the \emph{marked length
spectrum} $\Lambda_\alpha$
to be the set of lengths of all periodic
geodesics representing $\alpha$.
This set is closed and nowhere dense in $\R$.
Hence $\ell_\alpha:=\inf\Lambda_\alpha$
belongs to $\Lambda_\alpha$.
Let the open unit disc bundle
$DT^*M\subset T^*M$ be equipped
with the canonical
symplectic form $\omega_0$.

\begin{theorem}[\cite{JOA4}]
\label{thm:exist-orbit}
Let $\alpha$ be a homotopy class
of free loops in $M$.
Then every compactly supported
Hamiltonian
$H\in\Cinf([0,1]\times DT^*M)$
satisfying
$$
     \sup_{[0,1]\times M} H
     =:-c\le-\ell_\alpha
$$
admits a 1-periodic orbit $z=(x,y)$
with $[x]=\alpha$ and
action $\Aa_H(z)\ge c$.
\end{theorem}

\begin{figure}[ht]
\parbox[b]{3cm}{
  \epsfig{figure=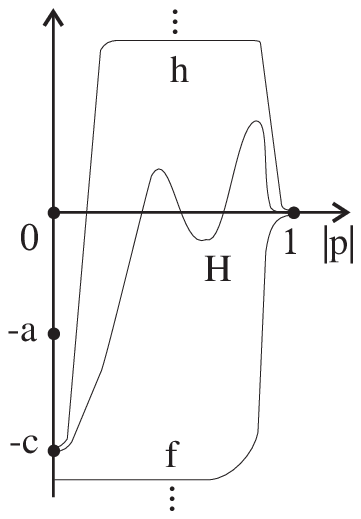,width=\linewidth}
}
\hspace{-.2cm}
\parbox[b]{9.1cm}{
  \begin{equation*}
  \xymatrix{
    &
     \underrightarrow{SH}^{(a,\infty);c;\alpha}_*
     \ar[r]^\simeq_{a\in(0,c]}
    &
     H_*(\Ll_\alpha M)
    \\
     HF_*^{(a,\infty)}(H;\alpha)
     \ar[ur]^{\iota_H} 
     \ar @{} [r] |{\qquad\circlearrowleft}
    &
     \ar @{} [r] |{\circlearrowleft}
    &
    \\
    &
    \underleftarrow{SH}^{(a,\infty);\alpha}_*
     \ar[ul]^{\pi_H}
     \ar[uu]_T
     \ar[r]^\simeq_{a\in\R^+\setminus
                    \Lambda_\alpha}
                  &
     H_*(\Ll^{\frac{a^2}{2}}_\alpha M)
     \ar[uu]_{\overset{\scriptstyle
              [\iota]\not=0\quad\;\;}
              {\text{if $a\ge\ell_\alpha$}}}
  }
  \end{equation*}
}
\caption{The Hamiltonian sandwich.}
\label{fig:fig-idea}
\end{figure}

The idea of proof is to sandwich
the Hamiltonian $H$ between two
Hamiltonians $f$ and $h$ (see
Figure~\ref{fig:fig-idea}) whose
action filtered
Floer homologies are computable and nonzero.
Then prove that the
\emph{monotonicity homomorphism}
associated to $f$ and $h$
is nonzero. This is a homomorphism
from Floer homology of a given Hamiltonian
to Floer homology of any pointwise larger
Hamiltonian. Its crucial properties are,
firstly, it respects the action window
and, secondly, it factors through Floer
homology of any third intermediate Hamiltonian.
In our case it factors
through $\HF_*^{(a,\infty)}(H;\alpha)$,
which therefore
is nonzero. Roughly speaking,
if $c$ is a regular
value of $\Aa_H|_{\Ll_\alpha T^*M}$,
set $a:=c$ and we are done.

The key tool to calculate
the Floer homologies of $f$ and $h$
is a refinement of the
Floer continuation principle.
Namely, two Hamiltonians connected
by an \emph{action regular homotopy}
(meaning that the boundary of the action
window consists of regular values
throughout the homotopy)
have the same Floer homology.
Action regularity is most easily
checked for radial Hamiltonians.
\begin{remark}[Radial  Hamiltonians]
\label{rmk:radial}\rm
Let $h:\R\to\R$
be a smooth symmetric function.
A Hamiltonian of the form
$H(t,q,p)=h(\abs{p})$
is called \emph{radial}.
Whenever the slope of $h$ at
a point $r_\ell$
is equal to an element
$\ell\in\Lambda_\alpha$,
then there is a 1-periodic orbit $z_\ell$
of the Hamiltonian flow 
on the sphere bundle of radius
$r_\ell$ and all 1-periodic
orbits arise that way.
Moreover, the symplectic
action of $z_\ell$ equals
minus the intercept of the tangent
to the graph at the point $r_\ell$.
\end{remark}
Hence we choose the Hamiltonians
$f$ and $h$ radial.
Then in each case we construct an
action regular homotopy towards
a \emph{convex} radial Hamiltonian,
since for these
Floer homology is computed
in~\cite{JOA4} on the basis
of Theorem~\ref{thm:main}.

Of course, for different $H$
we may have to take new choices for
$f$ and $h$. A convenient book 
keeping tool
to deal with this problem
is to take the inverse limit
over all $f$ and the direct
limit over all $h$ subject to
the restriction
$\sup_{[0,1]\times M} h\le-c$.
These limits are called symplectic and
relative symplectic homology
and were introduced in~\cite{FH94,CFH95}
and~\cite{BPS03}, respectively.
The monotone homomorphisms
descend to a natural homomorphism $T$.
The main part of~\cite{JOA4} is devoted
to establish and prove commutativity
of the rectangular
part of the diagram in
Figure~\ref{fig:fig-idea}.
Here the homomorphism $[\iota]$
induced by inclusion
does not vanish
whenever $a\ge\ell_\alpha$.
The monotone homomorphisms
descend to a natural homomorphism
$T$. It follows that $T$
is nonzero whenever
$a\in[\ell_\alpha,c]\setminus\Lambda_\alpha$.
In the case that $c$
is a regular value of
$\Aa_H|_{\Ll_\alpha T^*M}$
and $-c<-\ell_\alpha$
we can choose $a=c$
and are done.
Otherwise choose a sequence
of Hamiltonians $H_\nu$
converging to $H$ in $C^\infty$ and
such that the corresponding $c_\nu$
satisfy both requirements above.
From the resulting
sequence of periodic orbits
extract a subsequence
whose limit is the periodic orbit
claimed by Theorem~\ref{thm:exist-orbit}.

Consider the radial Hamiltonian
whose graph in $\R^2$ consists
of a straight line from
$(0,-\ell_\alpha)$ to $(1,0)$
and which is zero elsewhere.
Approximate it by a smooth
function all of whose slopes are
strictly less then $\ell_\alpha$
and which therefore does not
admit any 1-periodic
orbit representing $\alpha$.
This example shows that the condition
in Theorem~\ref{thm:exist-orbit}
is sharp. It follows that the relative
capacity defined in~\cite{BPS03}
and associated to
$(DT^*M,M,\alpha)$
equals $\ell_\alpha$.
As a byproduct we obtain
in~\cite{JOA4}
a multiplicity version of
the Weinstein conjecture
for compact hypersurfaces
$Q\subset T^*M$ of contact type
which enclose
the zero section $M$.
More precisely, for every nontrivial
$\alpha$ we obtain existence
of a closed characteristic
on $Q$ whose projection to $M$
represents $\alpha$.
In the nonsimply connected case
this refines the result of Hofer and
Viterbo~\cite{HV88}
in the case $\alpha=0$.
Viterbo informed us that
their techniques should also provide
multiplicities.



\begin{thebibliography}{9999}
\small

\bibitem[AM04]{AM04}  A.~Abbondandolo and P.~Majer,
      Lectures on the Morse complex
      for infinite dimensional manifolds,
      Summer School on
      {\it Morse theoretic methods in
      non-linear analysis and symplectic
      topology}, Kluwert, Montreal~2004.

\bibitem[AS04]{AS04}  A.~Abbondandolo and M.~Schwarz,
      On the Floer homology of cotangent bundles,
      Preprint, August~2004.
      arxiv:math.SG/0408280, to appear in
      {\it Comm. Pure Appl. Math.}

\bibitem[AS04b]{AS04b}  A.~Abbondandolo and M.~Schwarz,
      Floer homology of cotangent bundles,
      Summer School on
      {\it Morse theoretic methods in
      non-linear analysis and symplectic
      topology}, Kluwert, Montreal~2004.

\bibitem[BPS03]{BPS03}  P.~Biran, L.~Polterovich
      and D.A.~Salamon,
      Propagation in Hamiltonian dynamics
      and relative symplectic homology,
      {\it Duke Math. J.}
      {\bf 119} (2003), 65--118.

\bibitem[Ch84]{Ch84}  M.~Chaperon,
      Une id\'ee du type
      ``g\'eod\'esiques bris\'ees''
      pour les syst\`emes hamiltoniens,
      {\it C. R. Acad. Sc. Paris}
      {\bf 298} (1984), 293--96.

\bibitem[CS99]{CS99}  M.~Chas and D.~Sullivan,
      String topology,
      Preprint~1999.
      arxiv:math.GT/ 9911159, to appear in
      {\it Annals of Mathematics}.

\bibitem[C94]{C94}  K.~Cieliebak,
      Pseudo-holomorphic curves and
      periodic orbits on cotangent bundles,
      {\it J. Math. Pures Appl.} {\bf 73}
      (1994), 251--78.

\bibitem[C02]{C02}  K.~Cieliebak,
      Handle attaching in symplectic
      homology and the chord conjecture,
      {\it J. Eur. Math. Soc.} {\bf 4}
      (2002), 115--42.

\bibitem[Co04]{Co04}  R.~Cohen,
      Lectures on Morse theory, graphs,
      and string topology, Summer School on
      {\it Morse theoretic methods in
      non-linear analysis and symplectic
      topology}, Kluwert, Montreal~2004.
      arxiv:math.GT/0411272

\bibitem[CFH95]{CFH95}  K.~Cieliebak, 
      A.~Floer and H.~Hofer,
      Symplectic homology, II.
      A general construction,
      {\it Math.~Z.} {\bf 218} (1995),
      103--22.

\bibitem[F89]{F89}  A.~Floer,
      Symplectic fixed points and
      holomorphic spheres,
      {\it Comm. Math. Phys.} {\bf 120} (1989),
      575--611.

\bibitem[F89b]{F89b}  A.~Floer,
      Witten's complex and infinite
      dimensional Morse theory,
      {\it J. Differential Geom.}
      {\bf 30} (1989), 207--21.

\bibitem[FH94]{FH94}  A.~Floer and H.~Hofer,
      Symplectic homology,
      I. Open sets in $\C^n$,
      {\it Math.~Z.} {\bf 215} (1994),
      37--88.

 \bibitem[GL00]{GL00}  D.~Gatien and F.~Lalonde,
      Holomorphic cylinders with
      Lagrangian boundaries and
      Hamiltonian dynamics,
      {\it Duke Math. J.}
      {\bf 102} (2000), 485--511.

\bibitem[G89]{G89}  A.B.~Givental,
      Periodic maps in symplectic topology,
      {\it Funct. Anal. Appl.}
      (English translation)
      {\bf 23} (1989), 287--300.

\bibitem[HV88]{HV88}  H.~Hofer and C.~Viterbo,
      The Weinstein conjecture in cotangent
      bundles and related results,
      {\it Annali Sc. Norm. Sup. Pisa},
      Serie IV, Fasc. III,
      {\bf 15} (1988), 411--45.

\bibitem[L04]{L04}  F.~Laudenbach,
      Symplectic geometry and
      Floer homology, In {\it
      Symplectic geometry and
      Floer homology.
      A survey of Floer homology for
      manifolds with contact type
      boundary or symplectic homology} 1--50, 
      Ensaios Mat. {\bf 7},
      Soc. Brasil. Mat.,
      Rio de Janeiro, 2004.

\bibitem[L03]{L03}  Y.-J.~Lee,
      Non-contractible periodic orbits,
      Gromov invariants, and
      Floer-theoretic torsions,
      Preprint 2003. math.SG/0308185

\bibitem[MS04]{MS04}  D.~McDuff and D.A.~Salamon,
      {\it $J$-holomorphic curves and Symplectic Topology},
      Colloquium Publications, Vol. {\bf 52},
      American Mathematical Society, Providence,
      Rhode Island, 2004.

\bibitem[M64]{M64}  J.~Milnor,
      {\it Morse theory},
      Annals of Mathematics Studies 51,
      Princeton University Press,
      Princeton NJ 1964.

\bibitem[O04]{O04}  A.~Oancea,
      A survey of Floer homology for
      manifolds with contact type
      boundary or symplectic homology,
      In {\it Symplectic geometry and
      Floer homology.
      A survey of Floer homology for
      manifolds with contact type
      boundary or symplectic homology} 51--91,
      Ensaios Mat. {\bf 7},
      Soc. Brasil. Mat.,
      Rio de Janeiro, 2004.

\bibitem[S97]{S97}  D.A.~Salamon,
      Lectures on Floer Homology,
      In {\it Symplectic Geometry
      and Topology} 143--230,
      edited by Y.~Eliashberg and L.~Traynor,
      IAS/Park City Mathematics Series,
      Vol {\bf 7}, 1999.

\bibitem[SW03]{JOA3}  D.A.~Salamon and J.~Weber,
      Floer homology and the heat flow,
      Preprint, ETHZ, April 2003,
      revised September~2004.
      arxiv:math.SG/0304383, to appear in
      {\it Geom. Funct. Anal.}

\bibitem[Sch93]{Sch93}  M.~Schwarz,
      {\it Morse homology}, PM {\bf 111},
      Birkh\"auser, Basel 1993.

\bibitem[T94]{T94}  L.~Traynor,
      Symplectic homology
      via generating functions,
      {\it Geom. Funct. Anal.}
      {\bf 4} (1994), 718--48.

\bibitem[V95]{V95}  C.~Viterbo,
      Generating functions in
      symplectic topology and applications,
      In {\it Proceedings ICM~94 Z\"urich}
      537--47, Vol~{\bf 1}, Birkh\"auser,
      Basel~1995.

\bibitem[V96]{V96}  C.~Viterbo,
      Functors and computations
      in Floer homology
      with applications, Part II,
      Preprint October~1996, 
      revised~2003. \\
      http://math.polytechnique.fr/cmat/viterbo/viterbo.html

\bibitem[V97]{V97}  C.~Viterbo,
      Exact Lagrange submanifolds,
      periodic orbits and the
      cohomology of free loop spaces,
      {\it J. Differential Geom.}
      {\bf 47} (1997), 420--68.

\bibitem[V99]{V99}  C.~Viterbo,
      Functors and computations in
      Floer homology with
      applications, I,
      {\it Geom. Funct. Anal.}
      {\bf 9} (1999), 985--1033.

\bibitem[W96]{JOA0c} J.~Weber,
      Morse theory on the loop space of
      flat tori and symplectic Floer theory,
      SFB~288 {\sl Geometrie \& Quantenphysik},
      Preprint~237, TU Berlin 1996.
      arxiv:math.dg-ga/9612012

\bibitem[W99]{JOA1}  J.~Weber,
      {\it $J$-holomorphic curves in
      cotangent bundles and the heat flow},
      PhD thesis, TU Berlin, 1999. 

\bibitem[W02]{JOA2}  J.~Weber,
      Perturbed closed geodesics are periodic
      orbits: Index and transversality, 
      {\it Math. Z.} {\bf 241}
      (2002), 45--81.

\bibitem[W04]{JOA4}  J.~Weber,
      Noncontractible periodic orbits in 
      cotangent bundles and Floer homology, 
      Preprint, ETHZ, April 2004.
      arxiv:math.SG/0410609, to appear in
      {\it Duke Math. J.}

\bibitem[W04b]{JOA5}  J.~Weber,
      The Morse-Witten complex via dynamical systems,
      Preprint, Universit\"at M\"unchen,
      November 2004.
      arxiv:math.GT/0411465, to appear in
      {\it Expo. Math.}

\bibitem[W]{JOA-FUTURE}  J.~Weber,
      The heat flow and the 
      homology of the loop space,
      In preparation.

\bibitem[Wi82]{Wi82}  E.~Witten,
      Supersymmetry and Morse theory,
      {\it J. Differential Geom.}
      {\bf 17} (1982), 661--92.
      
\end{thebibliography}
\end{document}